\input amstex
\documentstyle{amsppt}
\magnification=\magstep1
\hsize=5in
\vsize=7.3in
\TagsOnRight
\topmatter
\title Strong additivity and conformal nets
\endtitle
\author Feng  Xu \endauthor
\address{Department of Mathematics, University of California at Riverside,
2208 Sproul Hall, CA 92521}
\endaddress
\email{xufeng\@ math.ucr.edu}
\endemail
\abstract{     
We show that the fixed point subnet of a strongly additive conformal net
under the action of a compact group is strongly additive. Using the idea
of the proof we define the notion of strong additivity for a pair of
conformal nets and we show that a key result  about the induction of
the pair which we proved previously under the finite index assumption
can be generalized to strongly additive pairs of conformal nets.
These results are applied to classify  conformal nets of
central charge $c=1$ which are not necessarily rational and satisfy
a spectrum condition.}
\endabstract
\dedicatory{ Dedicated to Masamichi Takesaki on the occasion of his 70th
birthday}
\enddedicatory
\thanks   
I'd like to thank Professors Chongying Dong, 
J\"{u}rgen Fuchs, Christoph Schweigert
, Yasuyuki Kawahigashi for useful comments, and especially Professor
Roberto Longo for pointing out reference [LRo] which helps to
improve Prop. 3.4. 
This work is partially
supported by a NSF grant .
\endthanks
\endtopmatter
\document

\heading \S1. Introduction \endheading 
Let ${\Cal A}$ be a conformal net (or precosheaf) (cf. \S2.1).  
Let $I_1,I_2$ be two disjoint intervals of the circle. 
In this paper  an interval of the circle is defined to be 
an open connected
proper subset of the circle.  By [FJ] ${\Cal A}$ is {\it additive}
in the sense that if $I_n$ is a sequence of intervals which cover
an interval $J,$ i.e., 
$\cup_nI_n \supset J,$ then
$
\vee_n {\Cal A}(I_n) \supset  {\Cal A}(J)
$
where $\vee_n  {\Cal A}(I_n)$ is the von Neumann algebra generated by
${\Cal A}(I_n), n\in {\Bbb N}.$
${\Cal A}$ is {\it strongly additive} if ${\Cal A}(I_1)\vee
{\Cal A}(I_2)= {\Cal A}(I)$ where $I_1\cup I_2$ is obtained
by removing an interior point from an interval $I$. A conformal
net which violates strong additivity is given in [BS].
\par
Strong additivity seems to be a rather technical condition, but it
plays an important role in studying the representations of conformal
nets (cf. [KLM], [Xo]).\par
In Prop. 2.8 of [Xo] we proved that 
the fixed point subnet of ${\Cal B}\subset {\Cal A}$
(cf. \S3) of a strongly 
additive net ${\Cal A}$ under the action of a finite group is 
strongly additive using Galois correspondence for the actions of finite
groups on von Neumann algebras. This result is generalized in [Ls] to
the case when   ${\Cal B}\subset {\Cal A}$ has finite Jones index using 
ideas of transplanting conformal subnets. Both proofs depend on the finite
index condition and it is not immediately clear how to generalize them to
the case when the finite group is replaced by a compact but infinite group.
The first result Theorem 2.4 of this paper is to show that Prop. 2.8
of [Xo] is true when the finite group is replaced by a compact group
using the results of [Ha]. As in [Xo], the idea is to show that
$$
{\Cal A}(I_1)\vee
{\Cal B}(I_2)= {\Cal A}(I) \tag 1
$$ where $I_1, I_2,I$ are as above, and we use
a dilation argument similar to [Xo]. But the proof is different from 
[Xo] and [Ls], and the idea of the proof also gives a simpler proof
of the strong additivity results of [Xo] and [Ls] (cf. Remark after Th. 2.4)
\par 
It turns out that (1) can be used to generalize a key  result 
Th. 3.3 about inductions 
between  ${\Cal B}$  and ${\Cal A}$ obtained in [Xb]. 
Since the results of
[Xb] have many applications (cf. [BE1-3], [FRS], [PZ]), 
we define a pair of nets
${\Cal B}\subset {\Cal A}$ to be strongly additive if (1) is satisfied
(cf. \S3.2) and we study such pairs in \S3.2.  In view of the applications in 
\S4, we prove Theorem 3.8 in \S3.2 which is a generalization of 
Th. 3.3 in [Xb], and the proof is 
different from the original proof in [Xb]. We note that
a few other results in [Xb] also generalize to our current setting and
we plan to return to them in a future publication.\par
In \S4 we apply the results of \S2 and \S3 to classify conformal nets 
${\Cal A}$
with 
central charge $c=1$. The idea is very simple and we describe it here
roughly as follows. The $c=1$ Virasoro subnet  ${\Cal B}\subset {\Cal A}$ 
(cf. \S4.1) has a representation (the Vector representation)
with Jones index $4$ (cf. Lemma 4.1), and 
since ${\Cal B}$ is
strongly additive by Th. 2.4,  the induced endomorphism to 
${\Cal A}$ has Jones index $4$ by Prop. 3.4 and Lemma 4.1.  
Hence the principal graph of this
induced endomorphism is one of the ``A-D-E'' graphs 
listed in [GHJ], and in fact 
such subfactors are classified in [Po1]. 
Depending on the nature of the principal graph,
we can obtain enough information on ${\Cal A}$ to identify 
 ${\Cal A}$ as the ``A-D-E'' list given in \S4.1 under a spectrum
condition in \S4.2. The spectrum condition states  that a degenerate
representation of the Virasoro net other than the vacuum representation 
must appear in  ${\Cal A}$ if ${\Cal A}\neq {\Cal B}$. 
This condition is  true
for all known examples and 
we conjecture that the spectrum condition is
always true. We note that a ``A-D-E'' type classification for
$c<1$ has been given in [KL], and there are some similarities between
our results and that of [KL]. But there are notable differences since
$c=1$ Virasoro net is not rational, and the results of \S2 and \S3
play a crucial role in our proofs. We also note that we have tried to
give different proofs for each cases in \S4.2 since we expect that
the ideas of these proofs, as well as the general results of 
\S2 and \S3, will have applications beyond those described in \S4. \par 
The rest of the  paper is organized as follows. In \S2  after reviewing 
the concept of conformal nets, 
subnets and relevant notions we prove the
strong additivity result Theorem 2.4 mentioned at the beginning of the 
introduction. 
In \S3.1 we consider a pair of conformal nets and extend some of the
results of [BE1-3], [LR] and [Xb] to this setting. 
Prop. 3.1 is essentially due to [LR], except under additivity assumptions
we have stronger properties which was 
proved under finite index conditions in [Xb]. Lemma 3.2, 3.3, Prop. 3.4
, 3,5 were also proved in [Xb] under finite index conditions, and are
generalized with suitable modifications of the original proofs to
our current case when the index is not necessarily finite. 
In \S3.2 we define the notion of strongly additive pair
. Lemma 3.6 and Prop. 3.7 give us many examples of such pairs. 
In Theorem 3.8 we generalize a key result Th. 3.3 of  [Xb] to 
strongly additive pair . Cor. 3.9  is a  direct consequence of 
Th. 3.8 and results of \S3.1. We note that Theorem 3.8 
and Cor. 3.9   give us powerful
tools to determine the nature of induced endomorphisms as already shown in
[Xb] under the finite index assumption.\par
In \S4 we apply the results in \S2 and \S3 to pairs    
${\Cal B}\subset {\Cal A}$ where ${\Cal B}$ is the Virasoro net with central
charge $c=1$. 
After reviewing the known ``A-D-E'' list of such pairs in \S4.1, 
we first determine a distinguished class of irreducible
representations of ${\Cal B}$ in Lemma 4.1, corresponding to the
irreducible
representations of $SU(2)$. The principal graph of the induced endomorphism
of the vector representation is the affine A-D-E graph described in
Lemma 4.2. In Lemma 4.3 we classify the covariant representations of
a simple net associated with a Heisenberg group. In Lemma 4.4 we determine
the  principal graph of the induced endomorphism
of the vector representation for the known list, showing that all possible
graphs do appear for conformal nets. As a by-product of Lemma 4.4 we 
determine  the fusion rules of $c=1$ Virasoro net where one of the
representation is  degenerate in Prop. 4.5, which contains the result in
[R1]. \S4.2 is devoted to the proof of Th. 4.6.  
We treat the discrete cases \S4.2.1-2 first where
the reconstruction results of [DR1-2] are used. We then treat the general case
according to the A-D-E type of the principal graph of the induced
vector representation, under the spectrum assumption defined in \S4.2.
Finally we consider one application of Th. 4.6 to a concrete example
appeared in our previous study of cosets. \par
After this paper appeared in the web, we have been informed by 
Professor Roberto Longo  that the 
results in  \S4.2.1-2  have some overlap with unpublished results 
of S. Carpi (cf. [K] ).
\heading \S2.  Strong additivity of subnets \endheading
\subheading{\S2.1 Sectors}
Let $M$ be a properly infinite factor
and  $\text{\rm End}(M)$ the semigroup of
 unit preserving endomorphisms of $M$.  In this paper $M$ will always
be the unique hyperfinite $III_1$ factors.
Let $\text{\rm Sect}(M)$ denote the quotient of $\text{\rm End}(M)$ modulo
unitary equivalence in $M$. We  denote by $[\rho]$ the image of
$\rho \in \text{\rm End}(M)$ in  $\text{\rm Sect}(M)$.\par
 It follows from
\cite{L3} and \cite{L4} that $\text{\rm Sect}(M)$, with $M$ a properly
infinite  von Neumann algebra, is endowed
with a natural involution $\theta \rightarrow \bar \theta $  ;
moreover,  $\text{\rm Sect}(M)$ is
 a semiring.\par
Let $\epsilon$ be a normal
faithful conditional expectation
from 
$M$ to $\rho(M)$.  We define a number $d_\epsilon$ (possibly
$\infty$) by:
$$
d_\epsilon^{-2} :=\text{\rm Max} \{ \lambda \in [0, +\infty)|
\epsilon (m_+) \geq \lambda m_+, \forall m_+ \in M_+
\}$$ (cf. [PP]).\par
 We define
$$
d_\rho = \text{\rm Min}_\epsilon \{ d_\epsilon \}.
$$   $d_\rho$ is called the statistical dimension of  $\rho$. 
$d_\rho^2$ will be called the (minimal) index of  $\rho$.
It is clear
from the definition that  the statistical dimension  of  $\rho$ depends only
on the unitary equivalence classes  of  $\rho$.
The properties of the statistical dimension can be found in
[L1], [L3] and  [L4].\par
For $\lambda $, $\mu \in \text{\rm Sect}(M)$, let
$\text{\rm Hom}(\lambda , \mu )$ denote the space of intertwiners from
$\lambda $ to $\mu $, i.e. $a\in \text{\rm Hom}(\lambda , \mu )$ iff
$a \lambda (x) = \mu (x) a $ for any $x \in M$. When there are 
several von Neumann algebras, we will write 
$\text{\rm Hom}(\lambda , \mu )_M$ to indicate the dependence on $M$. 

$\text{\rm Hom}(\lambda , \mu )$  is a  vector
space and we use $\langle  \lambda , \mu \rangle$ to denote
the dimension of this space.  $\langle  \lambda , \mu \rangle$
depends
only on $[\lambda ]$ and $[\mu ]$. 
Moreover if $ \nu, \lambda$ and $ \mu$ has finite index, we have
$\langle \nu \lambda , \mu \rangle =
\langle \lambda , \bar \nu \mu \rangle $,
$\langle \nu \lambda , \mu \rangle
= \langle \nu , \mu \bar \lambda \rangle $ which follows from Frobenius
duality (See \cite{L2} ).  We will also use the following
notation: if $\mu $ is a subsector of $\lambda $, we will write as
$\mu \prec \lambda $  or $\lambda \succ \mu $.  A sector
is said to be irreducible if it has only one subsector. 
In this paper we will sometimes use $1$ to denote the identity
sector if there is no possible confusion.  
\par
If $\lambda$ is a sector with finite statistical dimension, 
the {\it principal graph} $\Gamma$ of $\lambda$ is
a bi-parti graph defined as follows. The even vertices of  $\Gamma$, 
denoted by $ \Gamma_0$, are labeled by 
the irreducible sectors of $(\bar\lambda\lambda)^n, n\in {\Bbb N}$, and
the odd vertices of  $ \Gamma_1$ of $\Gamma$ are labeled by the  
irreducible sectors of $(\bar\lambda\lambda)^n\lambda, n\in {\Bbb N}$. 
An even vertex $x$ is connected to an odd vertex $y$ by   
$\langle x\lambda,y \rangle $ edges. We say that  $\lambda$ has {\it finite
depth} if  $\Gamma$ is a finite graph. 
Following  [Po2] we say that  $\lambda$ is {\it amenable} if
$||\Gamma||= d_\lambda$, where $\Gamma$ is considered 
a linear map  from $l^2(\Gamma_0)$ to  $l^2(\Gamma_1)$ as in [Po2].  

\subheading{\S2.2 Conformal nets and subnets}
In this section we recall the notion of irreducible conformal net 
(precosheaf)
and its  representations as described in [GL1].\par
By an {\it interval} we shall always mean an open connected subset $I$
of $S^1$ such that $I$ and the interior $I' $ of its complement are
non-empty.  We shall denote by  ${\Cal I}$ the set of intervals in $S^1$.
We shall denote by  $PSL(2, {\bold R})$ the group of
  conformal transformations on the complex plane
that preserve the orientation and leave the unit circle $S^1$ globally
invariant.  Denote by ${\bold G}$
the universal covering group of $PSL(2, {\bold R})$.  Notice that  ${\bold G}$
is a simple Lie group and has a natural action on the  unit circle
$S^1$. \par 
Denote by  $R(\vartheta )$  the (lifting to ${\bold G}$ of the) rotation by
an angle $\vartheta $. This one-parameter subgroup of
${\bold G}$ will be referred to as rotation group (denoted by Rot)
in the following.
We may associate a
one-parameter group with any interval $I$ in the following way.
Let $I_1$ be the upper
semi-circle, i.e. the interval
$\{e^{i\vartheta }, \vartheta \in (0, \pi )\}$.
 By using the Cayley transform
$C:S^1 \rightarrow {\bold R} \cup \{\infty \}$ given by
$z\rightarrow -i(z-1)(z +1)^{-1}$,
we may identify  $I_1$
with the positive real line ${\bold R}_+$. Then we consider
the one-parameter group $\Lambda _{I_1}(s)$ of diffeomorphisms of
$S^1$  such that
$$
C\Lambda _{I_1} (s) C^{-1} x = e^s x \, ,
\quad  s, x\in {\bold R} \, . \,  
C\Lambda _{I_1} (s) C^{-1} x = e^s x \, ,
\quad  s, x\in {\bold R} \, .
$$
We also associate with $I_1$ the reflection $r_{I_1}$ given by
$
r_{I_1}z = \bar z
$
where $\bar z$ is the complex conjugate of $z$.  It follows from
the definition  that
$\Lambda _{I_1}$ restricts to an orientation preserving diffeomorphisms of
$I_1$, $r_{I_1}$ restricts to an orientation reversing diffeomorphism of
$I_1$ onto $I_1^\prime $. \par
Then, if $I$ is an interval and we choose $g\in {\bold G}$ such that
$I=gI_1$ we may set
$$
\Lambda _I = g\Lambda _{I_1}g^{-1}\, ,\qquad
r_I = gr_{I_1}g^{-1}\, .
$$      
Let $r$ be an orientation reversing isometry of $S^1$ with
$r^2 = 1$ (e.g. $r_{I_1}$).  The action of $r$ on $PSL(2, {\bold R})$ by
conjugation lifts to an action $\sigma _r$ on ${\bold G}$, therefore we
may consider the semidirect product of
${\bold G}\times _{\sigma _r}{\bold Z}_2$.   Since
${\bold G}\times _{\sigma _r}{\bold Z}_2$ is a covering of the group
generated by $PSL(2, {\bold R})$ and $r$,
${\bold G}\times _{\sigma _r}{\bold Z}_2$ acts on $S^1$. We call
(anti-)unitary a representation $U$ of
${\bold G}\times _{\sigma _r}{\bold Z}_2$ by operators on ${\Cal H}$ such
that $U(g)$ is unitary, resp. antiunitary, when $g$ is orientation
preserving, resp. orientation reversing. \par
Now we are ready to define a  conformal
net (precosheaf). \par
An irreducible
conformal net ${\Cal A}$ of von Neumann
algebras on the intervals of $S^1$ 
is a map
$$
I\rightarrow {\Cal A}(I)
$$
from ${\Cal I}$ to the von Neumann algebras on a separable Hilbert space
${\Cal H}$ that verifies the following properties:
\vskip .1in
\noindent
{\bf A. Isotony}.  If $I_1$, $I_2$ are intervals and
$I_1 \subset I_2$, then
$$
{\Cal A}(I_1) \subset {\Cal A}(I_2)\, .
$$
 
\vskip .1in
\noindent
{\bf B. Conformal invariance}.  There is a nontrivial unitary
representation $U$ of
${\bold G}$  on         
${\Cal H}$ such that
$$
U(g){\Cal A}(I)U(g)^* = {\Cal A}(gI)\, , \qquad
g\in {\bold G}, \quad I\in {\Cal I} \, .
$$
\vskip .1in
\noindent
{\bf C. Positivity of the energy}.  The generator of the rotation subgroup
$U(R(\vartheta ) )$ is positive.
 
\vskip .1in
\noindent
{\bf D.  Locality}.  If $I_0$, $I$ are disjoint intervals then
${\Cal A}(I_0)$ and $A(I)$ commute.
 
The lattice symbol $\vee $ will denote `the von Neumann algebra generated
by'.
\vskip .1in
\noindent                                               
{\bf E. Existence of the vacuum}.  There exists a unit vector
$\Omega $ (vacuum vector) which is $U({\bold G})$-invariant and cyclic for
$\vee _{I\in {\Cal I}}{\Cal A}(I)$.
\vskip .1in
\noindent
{\bf F. Irreducibility}.  The only
$U({\bold G})$-invariant vectors are the scalar multiples of $\Omega$.
 
\vskip .1in
\noindent
The term irreducibility is due to the fact (cf. Prop. 1.2 of [GL1]) that
under the assumption of {\bf F} $\vee_{I\in {\Cal I}} A(I) =B({\Cal H})$.
\par       
We have the following (cf. Prop. 1.1 of [GL1]):
\proclaim{2.1 Proposition} Let ${\Cal A}$ be an irreducible 
conformal net. The 
following hold:
\roster
\item"{(a)}"  Reeh-Schlieder theorem: $\Omega $ is cyclic and separating 
for each von Neumann algebra ${\Cal A}(I)$, $I\in {\Cal I}$.
\item"{(b)}"  Bisognano-Wichmann property: $U$ extends to an 
(anti-)unitary representation of ${\bold G}\times _{\sigma _r}{\bold Z}_2$ 
such that, for any $I\in {\Cal I}$,
$$
\align
U(\Lambda _I(2\pi t)) &= \Delta _I^{it} \, \\
U(r_I) &= J_I \, 
\endalign
$$
where $\Delta _I$, $J_I$ are the modular operator and the modular 
conjugation associated with $({\Cal A}(I), \Omega )$ .
For each $g\in {\bold G} \times _{\sigma _r} {\bold Z}_2$
$$
U(g){\Cal A}(I)U(g)^* = {\Cal A}(gI) \, .
$$
\endroster
\roster
\item"{(c)}"  Additivity: if a family of intervals $I_i$ covers the 
interval $I$, then
$$
{\Cal A}(I) \subset \vee _i {\Cal A}(I_i)\, .
$$
\item"{(d)}" Haag duality: ${\Cal A}(I)' = {\Cal A}(I')$
\endroster
\endproclaim
A  {\it representation} $\pi $ of
${\Cal A}$ is a family of representations $\pi _I$ of the
von Neumann algebras ${\Cal A}(I)$, $I\in {\Cal I}$, on a
separable Hilbert space ${\Cal H}_\pi $ 
such that 
$$
I\subset \bar I \Rightarrow \pi _{\bar I} \mid _{{\Cal A}(I)}
= \pi _I \quad \text{\rm (isotony)} 
$$
A {\it covariant  representation} $\pi $ of
${\Cal A}$ is a family of representations $\pi _I$ of the
von Neumann algebras ${\Cal A}(I)$, $I\in {\Cal I}$, on a
separable Hilbert space ${\Cal H}_\pi $ and a unitary representation
$U_\pi $ of the covering group ${\bold G}$ of $PSL(2, {\bold R})$
with positive energy such that the following properties hold:
$$
\align
I\subset \bar I \Rightarrow \pi _{\bar I} \mid _{{\Cal A}(I)}
= \pi _I \quad &\text{\rm (isotony)} \\
\text{\rm ad}U_\pi (g) \cdot \pi _I = \pi _{gI}\cdot
\text{\rm ad}U(g) &\text{\rm (covariance)}\, .
\endalign
$$
A covariant representation $\pi$ is called irreducible if
 $\vee _{I\in {\Cal I}}\pi({\Cal A}(I)) = B({\Cal H}_\pi)$. By our definition
the irreducible conformal net
is in fact an irreducible representation  
of itself
and we will call this representation the {\it vacuum representation}. \par
We note that by [GL2] if a representation has finite index, then it 
is covariant. \par
Let $H$ be a  connected simply-laced compact Lie group. By Th. 3.2
of [FG], 
the vacuum positive energy representation of the loop group
$LH$ (cf. [PS]) at level $k$ 
gives rise to an irreducible conformal net 
denoted by {\it ${\Cal A}_{H_k}$}. By Th. 3.3 of [FG], every 
irreducible positive energy representation of the loop group
$LH$ at level $k$ gives rise to  an irreducible covariant representation 
of ${\Cal A}_{H_k}$.    
We also note that the vacuum representation
of the Virasoro algebra with central charge $c_0>0$ also give rise to
a conformal net denoted by ${\Cal A}_{c=c_0}$ (cf. \S3 of [FG]).
We will see such examples in \S4.\par 
Next we   recall some definitions from [KLM] .
As in [GL1] by an interval of the circle we mean an open connected
proper subset of the circle. If $I$ is such an interval then
$I'$ will denote the interior of the complement of $I$ in the circle.
We will denote by ${\Cal I}$ the set of such intervals.
Let $I_1, I_2\in {\Cal I}$. We say that $I_1, I_2$ are disjoint if
$\bar I_1\cap \bar I_2=\emptyset$, where $\bar I$
is the closure of $I$ in $S^1$.  
When $I_1, I_2$ are disjoint, $I_1\cup I_2$
is called a 1-disconnected interval in [Xj].  
Denote by ${\Cal I}_2$ the set of unions of disjoint 2 elements
in ${\Cal I}$. Let ${\Cal A}$ be an irreducible conformal net
as in \S2.1. For $E=I_1\cup I_2\in{\Cal I}_2$, let
$I_3\cup I_4$ be the interior of the complement of $I_1\cup I_2$ in 
$S^1$ where $I_3, I_4$ are disjoint intervals. 
Let 
$$
{\Cal A}(E):= A(I_1)\vee A(I_2), 
\hat {\Cal A}(E):= (A(I_3)\vee A(I_4))'.
$$ Note that ${\Cal A}(E) \subset \hat {\Cal A}(E)$.
Recall that a net ${\Cal A}$ is {\it split} if ${\Cal A}(I_1)\vee
{\Cal A}(I_2)$ is naturally isomorphic to the tensor product of
von Neumann algebras ${\Cal A}(I_1)\otimes
{\Cal A}(I_2)$ for any disjoint intervals $I_1, I_2\in {\Cal I}$.
${\Cal A}$ is {\it strongly additive} if ${\Cal A}(I_1)\vee
{\Cal A}(I_2)= {\Cal A}(I)$ where $I_1\cup I_2$ is obtained
by removing an interior point from $I$.
\proclaim{Definition 2.2 (Absolute rationality of [KLM])}
${\Cal A}$ is said to be absolute rational, or $\mu$-rational, if
${\Cal A}$ is split, strongly additive, and 
the index $[\hat {\Cal A}(E): {\Cal A}(E)]$ is finite for some
$E\in {\Cal I}_2$ . The value of the index
$[\hat {\Cal A}(E): {\Cal A}(E)]$ (it is independent of 
$E$ by Prop. 5 of [KLM]) is denoted by $\mu_{{\Cal A}}$
and is called the $\mu$-index of ${\Cal A}$. 
\endproclaim
Let  ${\Cal A}$ be a strongly additive  conformal net. 
This net is not directed. So when we discuss   
covariant representations
of ${\Cal A}$, we have to specify the intervals. To simply our notations, 
we fix a point $\xi \in S^1$. Denote by  ${\Cal I_\xi}$  the set of 
intervals $I$ with $\bar I\subset S^1-{\xi}$. Note that ${\Cal I_\xi}$
is a directed set under inclusion.  Let ${\Cal U}_{{\Cal A}}$
be the associated quasi-local $C^*$-algebra ${\Cal U}_{{\Cal A}}=\overline{
\cup_{J\in {\Cal I_\xi}} {\Cal A}(J)}$(norm closure).
We note that any representation of 
$\lambda$ of ${\Cal A}$ localized on $I$ restricts to a DHR endomorphism of
${\Cal U}_{{\Cal A}}$ localized on $I$ also denoted by $\lambda$ 
and vice versa 
(cf. Prop. 11 of [Ls]).  We will use these two descriptions interchangeably 
without further specifications. 
\par
Fix an $I\subset {\Cal I_\xi}$.
Let $\lambda,\mu$ be  representations of ${\Cal A}$. Choose
$I_-, I_+ \in {\Cal I_\xi}$ so that
$I_- \cap I= \emptyset= I_+ \cap I$, $I_-$ lies anti-clockwise to
$I$, and $I_+$ lies clockwise to
$I$. Choose  $\hat\lambda_+, \hat\lambda_-$ covariant representations
of ${\Cal A}$ unitarily equivalent to $\lambda$ but localized on
$I_+, I_-$ respectively and let $u_+, u_-$ be the unitary intertwinners.
Note that we will not distinguish local and global intertwinners since
they are the same when   ${\Cal A}$ is strongly additive.
The {\it braiding operators} are defined by
$$
\epsilon(\lambda,\mu):= \mu(u_+^*)u_+, 
\tilde\epsilon(\lambda,\mu)= \mu(u_-^*)u_-
$$
The properties of these operators are well known (cf. [Xb], [BE1]).
We note that $\epsilon(\lambda,\mu), \tilde\epsilon(\lambda,\mu)$ are
elements of $ {\Cal A}(I)$ and they are independent of the choices of
$u_+,u_-, I_+, I_-$ as long as 
$I_- \cap I= \emptyset= I_+ \cap I$, $I_-$ lies anti-clockwise to
$I$, and $I_+$ lies clockwise to $I$.  These operators satisfy Yang-Baxter
equation (YBE) and Braiding-Fusion equation (BFE), and we refer the
reader to [Xb] and [BE1] for more details.  \par

By a {\it conformal subnet} (cf. [Ls]) we shall mean a map
$$
I\in {\Cal I} \rightarrow {\Cal B}(I) \subset  {\Cal A}(I) 
$$ 
which associates to each interval $I\in {\Cal I}$ a von Neumann subalgebra
$ {\Cal B}(I)$ which is isotone
$$ {\Cal B}(I_1) \subset {\Cal B}(I_2), I_1\subset I_2,$$
and covariant with respect to the representation $U$, i.e.,
$$
U(g){\Cal B}(I)U(g)^* = {\Cal B}(gI)\, , \qquad
g\in {\bold G}, \quad I\in {\Cal I} \, .
$$
Note that when restricting to the intervals from the set ${\Cal I_\xi}$, 
the conformal nets 
${\Cal B}\subset{\Cal A}$ is a standard net of inclusions as defined in
\S3.1 of [LR]. We denote by $\gamma_{\Cal A}$ be the canonical 
endomorphism from  ${\Cal A}(J)$ to  ${\Cal B}(J), 
J\supset I, J\in {\Cal I_\xi }$ which is an   extension of canonical 
endomorphism from ${\Cal A}(I)$ to  ${\Cal B}(I)$ as defined by 
Cor. 3.3 of [LR]. 
The restriction of   $\gamma_{\Cal A}$ to
 ${\Cal B}$ will be simply denoted by $\gamma$, and when no confusion
arises, we will also denote  $\gamma_{\Cal A}$ by   $\gamma$. \par  
We have (cf. [KLM] or Prop. 2.4 of [Xo]):
\proclaim{Proposition 2.2}
Suppose ${\Cal B}\subset {\Cal A}$ is a standard net of inclusions
as defined in 3.1 of [LR]. Let $E\in {\Cal I}_2$.
If  ${\Cal A}\subset {\Cal B}$ has finite index denoted by 
$[{\Cal A}:{\Cal B}]$ and ${\Cal A}$ and ${\Cal B}$ are split, 
then
$$
[\hat {\Cal B}(E):{\Cal B}(E)] = [{\Cal A}:{\Cal B}]^2 
[\hat {\Cal A}(E):{\Cal A}(E)].
$$ 
\endproclaim

\subheading{\S2.3 Orbifolds}
Let ${\Cal A}$ be an irreducible conformal net on a Hilbert space
${\Cal H}$ and let $G$ be a compact group. Let $V:G\rightarrow U({\Cal H})$
be a faithful\footnotemark\footnotetext{
If $V:G\rightarrow U({\Cal H})$ is not faithful, we can take 
$G':= G/ker V$ and consider $G'$ instead.} 
unitary representation of $G$ on ${\Cal H}$.
\proclaim{ Definition 2.1} 
We say that $G$ acts properly on ${\Cal A}$ if the following conditions
are satisfied:\par
(1) For each fixed interval $I$ and each $g\in G$, 
$\alpha_g (a):=V(g)aV(g^*) \in {\Cal A}(I), \forall a\in
{\Cal A}(I)$; \par
(2) For each  $g\in G$ $V(g)\Omega = \Omega, \forall g\in G$.\par
\endproclaim
Suppose that a finite group $G$ acts properly on   ${\Cal A}$ as above.
For each interval $I$, define ${\Cal B}(I):=\{a\in {\Cal A}(I)| 
V(g)aV(g^*)=a, \forall g\in G \}$. Let $ {\Cal H}_0= \{x\in {\Cal H}|
V(g) x=x, \forall g\in G \}$ and $P_0$ the projection from ${\Cal H}$ to 
${\Cal H}_0$. Notice that $P_0$ commutes with every element of
${\Cal B}(I)$ and $U(g), \forall g\in {\bold G}$. 

Define ${\Cal A}^G(I):={\Cal B}(I)P_0$ on ${\Cal H}_0$. The unitary
representation $U$ of ${\bold G}$ on ${\Cal H}$ restricts to
an  unitary
representation (still denoted by $U$) of ${\bold G}$ on ${\Cal H}_0$.
By Prop. 2.2 of [Xo] 
the map $I\in {\Cal I}\rightarrow {\Cal A}^G(I)$ on $ {\Cal H}_0$ 
together with the  unitary
representation (still denoted by $U$) of ${\bold G}$ on ${\Cal H}_0$
is an
irreducible conformal net,
denoted by
${\Cal A}^G$ and will be called the {\it orbifold of ${\Cal A}$}
with respect to $G$. \par
The net ${\Cal B}\subset {\Cal A}$ 
is a standard net of inclusions (cf. [LR]) when restricting to 
intervals in ${\Cal I_\xi}$
with conditional expectation $E$ defined by 
$$
E(a):=\int_G  \alpha_g(a)dg, \forall a\in {\Cal A}(I)
$$ 
where $dg$ is the normalized Haar measure on $G$.
\proclaim{Lemma 2.3}
(1) For any interval $I$,
${\Cal A}^G(I)' \cap{\Cal A}(I) = {\Bbb C}; $\par
(2) Let $I$ be an interval, and $I_1,I_2$ are the connected components
of a set obtained from $I$ by removing an interior point of $I$. 
Let $g_n\in {\bold G}$ be a sequence of elements such that
$g_n I_1 = I_1$ and $g_n I_2$ is an increasing sequence intervals 
containing $I_2$, i.e.,
$I_2\subset g_n I_2\subset g_{n+1} I_2$, and
$\cup_n  g_n I_2=I_1'$ (One may take $g_n$ to be a sequence of
dilations). Let $x\in {\Cal B}(I_1)'\cap {\Cal A}(I_2')$, and suppose 
$y$ is a weak limit of a subsequence of
$Ad_{g_k}(x):= g_k x g_k^*$. Then $y= \langle x \Omega,\Omega\rangle id.$
\endproclaim
\demo{Proof:}
The first part is Prop. 2.1 of  [C] (also cf. [ALR] for related results). 
To prove the
second part, note that $Ad_{g_k}(x) \in {\Cal B}(I_1)'\cap {\Cal A}(g_kI_2)',$
and so if $y$ is a weak limit, 
$y\in  {\Cal B}(I_1)'\cap {\Cal A}(g_kI_2)',\forall k$, and it follows that
$y\in  {\Cal B}(I_1)'\cap {\Cal A}(I_1)'= {\Bbb C}$ by 
(1). On the other hand since $g_k\Omega=\Omega$, we have
$$
\langle y\Omega, \Omega\rangle
= \langle x\Omega, \Omega\rangle,
$$ and so $y= \langle x \Omega,\Omega\rangle id.$
\enddemo
\qed
\par
\proclaim{Theorem 2.4}
Let ${\Cal A}$ be an irreducible conformal net and let $G$
be a compact group acting properly on ${\Cal A}$. Suppose that 
${\Cal A}$ is  strongly additive.  Then ${\Cal A}^G$ is
also  strongly additive.
\endproclaim
\demo{Proof:}
Let $I$ be an interval, and $I_1,I_2$ are the connected components
of a set obtained from $I$ by removing an interior point of $I$. 
To show ${\Cal A}^G$ is strongly additive, it is sufficient to show that
${\Cal B}(I_1)\vee {\Cal B}(I_2) = {\Cal B}(I).$\par
Let us show that ${\Cal B}(I_1)\vee {\Cal A}(I_2)= {\Cal A}(I)$. 
Note that this will prove the theorem by simply applying $E$
to both sides of the equality.\par
Since
${\Cal A}$ is strongly additive, it is sufficient to show that
 ${\Cal B}(I_1)\vee {\Cal A}(I_2)= {\Cal A}(I_1)\vee {\Cal A}(I_2)$, or,
by taking commutants and using Haag duality
$$N:= {\Cal B}(I_1)'\cap {\Cal A}(I_2')= 
M:=  {\Cal A}(I_1')\cap {\Cal A}(I_2')
$$
Define $$ N_0:= {\Cal B}(I_1)' \supset M_0:={\Cal A}(I_1')$$
Notice that by Remark 4.5 of [I], $ N_0$ can be 
identified as the cross product
of $ M_0$ by $G$. By (a) of Th. 3.1 of [H], for each
continuous, positive definite function $\phi$ on $G$ there is a unique
$\sigma$-weakly continuous linear map $E_\phi$ on  $ N_0$
such that 
$$
E_\phi(m_0 x m_0') = m_0 E_\phi( x) m_0',
E_\phi(g)= \phi(g)g, \forall  m_0,  m_0'\in M_0, g\in G
$$ 
Let $g_n\in {\bold G}$ be a sequence of elements such that
$g_n I_1 = I_1$ and $g_n I_2$ is an increasing sequence intervals 
containing $I_2$, i.e.,
$I_2\subset g_n I_2\subset g_{n+1} I_2$, and
$\cup_n  g_n I_2=I_1'$ (One may take $g_n$ to be a sequence of
dilations).
Consider $$F_\phi(n_0):= Ad_{g_k^*} E_\phi( Ad_{g_k}n_0), \forall 
n_0\in N_0$$
Note that
$$
F_\phi(m_0 x m_0') = m_0 F_\phi( x) m_0',
F_\phi(g)= \phi(g)g, \forall  m_0,  m_0'\in M_0, g\in G
$$
Since $M_0 g M_0$ is weakly closed in $N_0$, it follows that
$
E_\phi=F_\phi,$ i.e.,
$$E_\phi( Ad_{g_k}n_0)=  Ad_{g_k}(E_\phi( n_0)), \forall 
n_0\in N_0$$
Let $x\in N$. Since ${\Cal A}(I_2)\subset {\Cal A}(I_1')$,
$E_\phi = \phi(1) id$ on ${\Cal A}(I_1')$, 
$E_\phi(x)\in N$.
Let $a,b\in M$. Then
$$
\align
\langle \frac{1}{\phi(1)} E_{\phi}(x)a\Omega, b\Omega \rangle
&= \langle \frac{1}{\phi(1)} E_{\phi}(b^*xa)\Omega, \Omega \rangle
\\
&=\langle \frac{1}{\phi(1)} Ad_{g_k}E_{\phi}(b^*xa)\Omega, \Omega \rangle
\\
= \langle \frac{1}{\phi(1)} 
E_{\phi}( Ad_{g_k}(b^*xa))\Omega, \Omega \rangle
\endalign
$$
Note that  by Lemma 2.3, there is a subsequence of
$ Ad_{g_k}b^*xa$ which converges weakly to 
$\langle b^*xa \Omega, \Omega\rangle$. Since $E_{\phi}$ is weakly continuous,
we must have
$$
\langle \frac{1}{\phi(1)} E_{\phi}(x)a\Omega, b\Omega \rangle
=\langle xa\Omega, b\Omega \rangle
$$
Since $M\Omega$ is dense in $H$, we have 
$E_{\phi}(x)= \phi(1)x.$
Let us choose $x$ to be a projection in $N$. By (b) of Th. 3.1 of
[H], 
$$
T(x)= sup_{\phi\ll \delta} E_{\phi}(x)
$$
where $T$ is the operator valued weights from $N_0$ to $M_0$, and
we write $\phi\ll \delta$ if $\phi$ is less than the Dirac measure in
the unit element of $G$ with respect to the ordering of positive definite
measures on $G$. 
So when $G$ is a finite group, 
$$
T(x)=x
$$ implying that $x\in M_0\cap N=M.$
When $G$ is an infinite group
$$
T(x)=x\cdot \infty
$$ in the notation of extended positive part of 
$M_0$ on P. 150 of [SS]. And so $x\in M_0\cap N=M.$
So we have shown that any projection of $N$ is a projection of
$M\subset N,$ and so $M=N$. 
\enddemo
\qed \par
We note that the above proof is different from the one given in [X0].
The same idea also gives a different proof of the result in \S3.5.2
of [Ls] under the assumption that 
${\Cal B}\subset{\Cal A}$ has finite index but 
without the assumption that ${\Cal A}$ is split if we modify the proof
as follows. Instead of using $E_\phi$ we use $E$ the minimal conditional
expectation from $N_0$ to $M_0$ which exists by finite index assumption.
Let us check that just like $E_\phi$, we have
$$E( Ad_{g_k}n_0)=  Ad_{g_k}(E( n_0)), \forall 
n_0\in N_0$$
Note that $Ad_{g_k^*}E_\phi( Ad_{g_k} \cdot)$ is a conditional expectation
from  $N_0$ to $M_0,$ and since $M_0\subset N_0$ is irreducible by 
Lemma 14 of [LS], we must have 
$Ad_{g_k^*}E( Ad_{g_k} \cdot)= E(\cdot).$ Now using
$E$ instead of $E_\phi$, the rest of the proof
goes through, and we get $E(x)=x, \forall x\in N,$ and so
$x\in N\cap M_0= M,$ i.e., $N=M$. \par
Let us consider a large class of  examples where Th. 2.4 can be applied. 
Let ${\Cal A}_{H_k}$
be the conformal net associated with representations of loop group
$LH$ at level $k$ as in \S2.2. Let $G\subset H$ be any closed subgroup.
By [TL]   ${\Cal A}_{H_k}$ is strongly additive, and it is easy
to check that $G$ or $G'$ as defined in the footnote of Definition 2.1
acts properly on  ${\Cal A}_{H_k}$. It follows that the fixed point net
${\Cal A}_{H_k}^{G'}$ is strongly additive by Th. 2.4. We will see
a special case of such examples in \S4.  
\heading \S3. Induction and strongly additive pairs\endheading 
\subheading{\S 3.1 Induction of a pair}
Let ${\Cal B}\subset {\Cal A}$ be a pair of conformal nets. In this
section we assume that ${\Cal B}$ is strongly additive. 
Fix a point $\xi \in S^1$. Denote by  ${\Cal I_\xi}$  the set of 
intervals $I$ with $\bar I\subset S^1-{\xi}$. Note that ${\Cal I_\xi}$
is a directed set under inclusion. Fix an $I\subset {\Cal I_\xi}$. 
All the intervals in this section will be in ${\Cal I_\xi}$ unless 
otherwise stated. Recall from \S2.2 that  $\gamma_{\Cal A}$ is the  
canonical endomorphism from  ${\Cal A}(J)$
into ${\Cal B}(J), J\supset I, J\in  {\Cal I_\xi}$ which is 
an extension of the canonical endomorphism from  ${\Cal A}(I)$
into ${\Cal B}(I)$  as given by Cor. 3.3 of [LR]. The restriction of
 $\gamma_{\Cal A}$ to  ${\Cal B}$ is denoted by $\gamma$.
We note that $\gamma$ is a (DHR) representation  of 
${\Cal B}$ which is unitarily equivalent to the defining representation
of ${\Cal B}$ on the vacuum Hilbert space of ${\Cal A}$ since
${\Cal B}$ is assumed to be strongly additive (cf. the proof of Prop. 17
in [Ls]).  
The following  is essentially Prop. 3.9 of [LR], except that under 
our conditions we have further properties.
\proclaim{Proposition 3.1}
Let ${\Cal B}\subset {\Cal A}$ be a  pair of conformal nets with 
 ${\Cal B}$ strongly additive. With every  
representation  $\lambda$ of ${\Cal B}$ associate
$$
\alpha_\lambda:= \gamma_{\Cal A}^{-1} Ad_\epsilon 
\lambda \gamma_{\Cal A}, 
\tilde\alpha_\lambda:= \gamma^{-1}_{\Cal A} Ad_{\tilde\epsilon} \lambda 
\gamma_{\Cal A} \tag 2
$$
where $\epsilon= \epsilon(\lambda, \gamma), \tilde\epsilon:= 
\tilde\epsilon(\lambda, \gamma)\lambda\gamma \rightarrow
\gamma \lambda$ are the braiding operators in 
${\Cal B}(I)$.  Then
$
\alpha_\lambda ({\Cal A}(I)) \subset 
{\Cal A}(I)
$
and $ \alpha_\lambda =\lambda $ on ${\Cal B}(I)$. 
Moreover,  $ \alpha_\lambda$ is localized on $I$ if and only if 
$$
\epsilon(\lambda, \gamma)  \epsilon(\gamma, \lambda)= id.
$$
\endproclaim
\demo{Proof} 
Choose a unitary intertwinner $u$ transporting $\lambda$ to $\hat\lambda$
localized on $I_+$ such that
$$
\epsilon(\lambda, \gamma)= \gamma(u^*)u
$$
Then for any $x\in {\Cal A}(I)
$ we have
$$
Ad_\epsilon \lambda\gamma(x) = Ad_{\gamma(u^*)} \hat\lambda\gamma =
\gamma( Ad_{u^*}x)
$$
Let $J$ be an interval in ${\Cal I_\xi}$ containing $I\cup I_+$. Then 
$\gamma( Ad_{u^*}x)\in \gamma( {\Cal A}(J
)),$
Note that the left hand side of the above depends on $J$ only through
the braiding operator 
$\epsilon$, and  by the invariance property of $\epsilon$ we can choose 
a decreasing sequence of intervals 
$J_n\supset \bar I$ such that $\cap_n J_n = I. $ By (c) of Prop. 2.1
$$
\cap_n {\Cal A}(J_n) = {\Cal A}(I).
$$ 
It follows that 
$$
\cap_n \gamma({\Cal A}(J_n)) = \gamma({\Cal A}(I)).
$$ 
So we have shown that
$$
Ad_\epsilon \lambda\gamma({\Cal A}(I) )\subset  \gamma({\Cal A}(I)).
$$
Hence 
$
\alpha_\lambda ({\Cal A}(I)) \subset 
{\Cal A}(I)
$ and 
$$
\alpha_\lambda(x)= u^*xu , \forall x\in {\Cal A}(I)\tag 3
$$
Note that a similar formula hold for $\tilde\alpha_\lambda$ if
we choose the unitary intertwinner accordingly. \par
The rest of the proof is the same as that of Prop. 3.9 of [LR].
\enddemo
\qed
\par
The notation $\alpha_\lambda$ was introduced in [BE1]. In [Xb], 
a slightly different induction $a_\lambda\in End({\Cal B}(I))$ was used 
motivated by certain questions in subfactors 
and the relations between these two are
given in [Xi]. 
Let us point out one basic relation.
Let $\rho\in End({\Cal B}(I))$ be  such that
$$\rho ({\Cal B}(I))= \gamma({\Cal A}(I)), \rho\bar\rho =\gamma, $$ then
$$
a_\lambda(\gamma(a))= \rho^{-1}(\gamma\alpha_\lambda(a))
, \forall a\in {\Cal A}(I)$$
We note that all the results of \S3 can be written in terms of $a_\lambda$
using the above relation.\par
The following lemma is implicitly contained in [LR]:
\proclaim{Lemma 3.2}
(1) If $x\in Hom(\lambda, \mu)_{\Cal B}$, then
$x\in Hom(\alpha_\lambda, \alpha_\mu)_{\Cal A}$; \par 
(2) $\alpha_{\lambda\mu} = \alpha_{\lambda}\alpha_\mu$; \par
(3) If $[\delta]= [\lambda]+ [\mu]$, then
$[\alpha_{\delta}] = [\alpha_{\lambda}]+[\alpha_\mu]$;\par
\endproclaim
\demo{Proof:}
Let $u_\lambda$ and $u_\mu$ be the unitary intertwinners from $\lambda$ and
$\mu$ localized on $I$ to  $\hat\lambda$ and
$\hat \mu$ localized on $I_+$ respectively. We can choose $I_+$ to share
only one boundary point with $I$. Since 
 $x\in Hom(\lambda, \mu)$, we have
$u_\mu x u_\lambda^*\in {\Cal B(I\cup I_+)}\cap  {\Cal B(I)}' = {\Cal B}(I_+)
$   
by strong additivity of $ {\Cal B}$. It follows that
$u_\mu x u_\lambda^*$ commutes with every element of $ {\Cal A}(I)$,
and (1) is proved. (2) and (3) follow from the definitions, 
YBE and BFE as in \S3 of [BE1].
\enddemo
\qed
\par
\proclaim{Lemma 3.3}
(1) If $\lambda$ is a representation of ${\Cal B}$ localized
on $I$ and $\sigma\in End({\Cal A}(I))$, then
$$
\langle \alpha_\lambda, \sigma \rangle_{{\Cal A}} \leq \langle 
\lambda, \gamma\sigma\rangle_{{\Cal B}}
$$ \par
(2) If $\alpha_\lambda, \alpha_\mu$ are localized on $I$, and denote by
$\epsilon(\alpha_\lambda, \alpha_\mu)$ the braiding operator 
of  $\alpha_\lambda, \alpha_\mu$
condsidered as (DHR) representations of ${\Cal A}$, then
$$
\epsilon(\alpha_\lambda, \alpha_\mu)=\epsilon(\lambda,\mu)
$$\par
(3) If ${\Cal B}\subset {\Cal C}\subset {\Cal A}$ 
such that ${\Cal C}\subset {\Cal A}$
and  ${\Cal B}\subset {\Cal C}$ are conformal subnets, and 
$\lambda_1:=\alpha^{{\Cal B}\rightarrow {\Cal C}}_\lambda$ is localized on
$I$. Then
$$
\alpha^{{\Cal C}\rightarrow {\Cal A}}_{\lambda_1}(a)  
= \alpha^{{\Cal B}\rightarrow {\Cal A}}_{\lambda}(a) , \forall a\in {\Cal A}(I) 
$$
where ${\Cal C}\rightarrow {\Cal A}, {\Cal B}\rightarrow {\Cal A}$ means the induction as defined in Prop. 3.1 from ${\Cal C}$ to ${\Cal A}$ and
from ${\Cal B}$ to ${\Cal A}$ respectively. 
\endproclaim
\demo{Proof:}
Ad(1):
Let $E: {\Cal A}(I)\rightarrow {\Cal B}(I)$ be the faithful conditional
expectation. Let $v\in Hom(id, \gamma)_{\Cal B}$ be the isometry such that
$
E(\cdot)= v^* \gamma(\cdot) v
$ (cf. [LR]). For any $x\in Hom(\alpha_\lambda, \sigma)_{\Cal A}$, 
it is easy to
check that $\gamma(x) v\in  Hom(\lambda, \gamma\sigma)_{\Cal B}$. 
To prove (1), we just have to show that the linear map
$$x\in  Hom(\alpha_\lambda, \sigma)_{\Cal A}\rightarrow \gamma(x) 
v\in  Hom(\lambda, \gamma\sigma)_{\Cal B}$$
is  one-to-one. Assume that $\gamma(x) v=0$, then
$$
E(x^*x)= v^* \gamma(x^*) v \gamma(x) v
=0
$$
It follows that $x=0$ since $E$ is faithful. 
\par
Ad (2) and (3):
Let $u_\lambda$ be a unitary intertwinner from $\lambda$ to $\hat\lambda$.
By (1) of Lemma 3.2 $u_\lambda$ is also a unitary intertwinner 
from $\alpha_\lambda$ to $\alpha_{\hat\lambda}$. So
$$
\epsilon(\alpha_\lambda, \alpha_\mu)= \alpha_\mu(u_\lambda)^* u_\lambda
=\mu(u_\lambda)^* u_\lambda= \epsilon(\lambda,\mu)
$$
and
$$
\alpha^{{\Cal C}\rightarrow {\Cal A}}_{\lambda_1}(a) 
=u_\lambda^* (a) u_\lambda
=\alpha^{{\Cal B}\rightarrow {\Cal A}}_{\lambda}(a), \forall a\in 
 {\Cal A}(I)
$$
by formula (3) in Prop. 3.1.
\enddemo
\qed
\par 
Note that by Prop. 3.1, $\alpha_\lambda\in End({\Cal A}(I))$. We will use
$d_{\alpha_\lambda}$ to denote the statistical dimension of $\alpha_\lambda$.
\proclaim{Proposition 3.4} If $\lambda$ has finite index, then:\par

(1) $d_\lambda= d_{\alpha_\lambda};$ \par

(2)  (Commuting Squares) 
Let $E$ be the conditional expectation from ${\Cal A}(I)$ to
${\Cal B}(I),$ and
 $F_\lambda$ the minimal conditional expectation from
$ {\Cal A}(I)\rightarrow \alpha_\lambda({\Cal A}(I))$ .  
Then $EF_\lambda=F_\lambda E$. 

\endproclaim
\demo{Proof:}
Ad (1): The proof is completed in three steps.\par
In the
first step we show that $d_\lambda\geq d_{\alpha_\lambda}$. 
We give two  different proofs of $d_\lambda\geq d_{\alpha_\lambda}$ . The first proof is longer but
use the additivity in a similar way as in the proof of Prop. 3.1. \par
Denote by $E_\lambda:  {\Cal B}(I)\rightarrow  \lambda_\epsilon({\Cal B}(I))$
the unique minimal conditional expectation. Let us first show that
$$
E_\lambda(\gamma ({\Cal A}(I)))\subset  \gamma ({\Cal A}(I)).
$$
Note that by [L1], there exists an isometry
$$
v\in Hom (id, \bar\lambda_{\bar\epsilon}\lambda_{\epsilon})
$$ 
such that 
$$E_\lambda(\cdot)= \lambda_{\epsilon}(v^*) \lambda_{\epsilon}
\bar\lambda_{\bar\epsilon}(\cdot)\lambda_{\epsilon}(v)  
$$
where 
$$
\epsilon:=\epsilon(\lambda,\gamma)= \gamma(u^*)u,
\bar\epsilon:=\epsilon(\bar\lambda,\gamma)= \gamma(\bar u^*)\bar u 
, \lambda_\epsilon:= Ad_\epsilon \lambda, \bar\lambda_{\bar\epsilon}:
= Ad_{\bar\epsilon} \bar\lambda
$$
and $u,\bar u\in {\Cal B}(J), J\supset I$ are unitary intertwinners as
in the definition of braiding operators (cf. \S2.2).  
we have 
$$
\lambda_{\epsilon}(b)= \gamma(u^*)b \gamma(u),
\lambda_{\epsilon}(b)= \gamma(\bar u^*) b\gamma(\bar u)  
$$
and so
$$
E_\lambda(\gamma(a))= \gamma(u^*) v^* \gamma(\bar u^* a\bar u)v \gamma(u) 
$$

Let us show that
$$
v \in \gamma ({\Cal A}(I)).
$$
Since $v\in Hom(id, \bar\lambda_{\bar\epsilon} \lambda_{\epsilon})$ we have
$$
vb=  \gamma(\bar u^* u^*) b\gamma(u \bar u) v, \forall b\in {\Cal B}(I)
$$   
and so $$ \gamma(u \bar u)v\in {\Cal B}(I)'\cap {\Cal B}(J)={\Cal B}(I_+)$$
where the equality follows from the strong additivity of ${\Cal B}$.
Hence $$v=  \gamma(\bar u^* u^*) \gamma(u \bar u)v \in 
\gamma(  {\Cal B}(J))$$ since $\gamma$ is localized on $I$. As in the proof 
of Prop. 3.1 we can choose a decreasing sequence of $J_n$ such that
$\cap_nJ_n=I$ and we have $v\in \cap_n\gamma({\Cal B}(J_n))=
\gamma({\Cal B}(I)).$
From the expressions for $E_\lambda(\gamma(a))$ above we have proved that
$$
E_\lambda(\gamma({\Cal A}(I)))\in \cap_n\gamma({\Cal A}(J_n))=
\gamma({\Cal A}(I))
$$
Let us show that
$$
\gamma({\Cal A}(I))\cap \lambda_\epsilon({\Cal B}(I))
=  \lambda_\epsilon\gamma({\Cal A}(I))
$$
If $ \lambda_\epsilon(b)= \gamma(a)$, then
$$
 \gamma( u^*) b  \gamma( u) = \gamma(a)
$$ and so
$$
b=  \gamma( ua  u^*)\in \gamma({\Cal A}(J))
$$
and by the same argument as above  
$$
b\in \gamma({\Cal A}(I))
$$
and this shows that
$$
\gamma({\Cal A}(I))\cap \lambda_\epsilon({\Cal B}(I))
=  \lambda_\epsilon\gamma({\Cal A}(I))
$$

Hence 
$$
E_\lambda (\gamma({\Cal A}(I)) \subset  \gamma({\Cal A}(I))\cap  
\lambda_\epsilon({\Cal B}(I))= \lambda_\epsilon\gamma({\Cal A}(I))
$$
So the minimal index of
$\lambda_\epsilon\gamma({\Cal A}(I))\subset  \gamma({\Cal A}(I))$
is less or equal to $d_\lambda^2$. Recall that
$
\lambda_\epsilon\gamma= \gamma(\alpha_\lambda)
$
and so the minimal index of
$\lambda_\epsilon\gamma({\Cal A}(I))\subset  \gamma({\Cal A}(I))$
is  $d_{\alpha_\lambda}^2$. So we have shown that
$d_{\alpha_\lambda}\leq d_\lambda. $ \par
Now we give a second proof of $d_{\alpha_\lambda}\leq d_\lambda$. 

Let $F_\lambda: {\Cal B}(I)\rightarrow \lambda({\Cal B}(I))$
be the minimal conditional expectation. 
Since $\lambda$ has finite index, there are two isometries 
$$R_{\lambda \bar\lambda}\in (id, \lambda \bar\lambda), 
R_{\bar\lambda \lambda}\in (id, \bar\lambda\lambda)
$$ in ${\Cal B}(I)$
such that
$$
R_{\lambda \bar\lambda}^*\lambda(R_{\bar\lambda \lambda}) 
=R_{ \bar\lambda\lambda}^*\bar\lambda(R_{\lambda \bar\lambda})
=\frac{1}{d(\lambda)}
$$

By [L1] (also cf. [LR]) $F_\lambda$ is induced by

an isometry $R_{\bar\lambda\lambda}\in Hom(id, \bar\lambda \lambda)$ such that
$$
F_\lambda(\cdot)=\lambda(R_{\bar\lambda\lambda})\lambda \bar\lambda(\cdot)
\lambda(R_{\bar\lambda\lambda}^*)
$$
By (1) of Lemma 3.2 we have 
$$
R_{\bar\lambda\lambda}\in Hom(id, \alpha_{\bar\lambda} \alpha_\lambda)
R_{\lambda \bar\lambda}\in Hom(id, \alpha_\lambda\alpha_{\bar\lambda})
$$
and it follows from the properties of $R_{\bar\lambda\lambda}, R_{\lambda\bar\lambda}$ that (cf. [LR])
$$
F_\lambda(\cdot)=\lambda(R_{\bar\lambda\lambda}) 
\alpha_\lambda\alpha_{\bar\lambda}(\cdot)
\lambda(R_{\bar\lambda\lambda}^*)
$$
is a conditional expectation from
$ {\Cal A}(I)\rightarrow \alpha_\lambda({\Cal A}(I))$
with index $d_\lambda^2$. It follows that $d_{\alpha_\lambda}\leq d_\lambda$.
\par
In the second step we show that if $\lambda$ is amenable (cf. \S2.1), then 
$d_\lambda=d_{\alpha_\lambda}$. \par 

Let $\Gamma$ be the principal graph of $\lambda$, and $\Gamma_0$ the
set of even vertices. 
By Lemma 3.2 and properties of statistical dimensions 
$V:=(d_{\alpha_x})_{x\in \Gamma_0}$
is a vector which verifies $\Gamma\Gamma^t V= d_{\alpha_\lambda} 
d_{\alpha_{\bar\lambda}} V$. \par
By Prop. 1.3.5 of [Po2],  $$
||\Gamma||^2 = \lim_{n\rightarrow \infty} 
((\Gamma\Gamma^t)^n \delta, \delta)^{\frac{1}{n}}
$$ where $\delta\in 
l^2(\Gamma_0)$ is a vector which is $1$ at the identity sector and
$0$ elsewhere. 
Hence
$$
\align
||\Gamma||^2 &= \lim_{n\rightarrow \infty} 
((\Gamma\Gamma^t)^n \delta, \delta)^{\frac{1}{n}} \\
& \leq \lim_{n\rightarrow \infty} ((\Gamma\Gamma^t)^n \delta, 
V)^{\frac{1}{n}}\\
&= \lim_{n\rightarrow \infty} (\delta, (\Gamma\Gamma^t)^n V)^{\frac{1}{n}}\\
&= \lim_{n\rightarrow \infty}
(\delta, (d_{\alpha_\lambda}d_{\alpha_{\bar\lambda}})^n V)^{\frac{1}{n}}\\
&=d_{\alpha_\lambda}d_{\alpha_{\bar\lambda}}
\endalign
$$
Since $\lambda$ is amenable, we have
$d_\lambda^2=||\Gamma||^2$, and it follows from above 
$d_\lambda^2\leq  d_{\alpha_\lambda}d_{\alpha_{\bar\lambda}}$.
By the first step  we must have $d_\lambda=  d_{\alpha_\lambda}$.
\par
Finally, by Theorem 5.31 and remarks on Page 122 of [LRo], 
$\lambda$ is always amenable. Combining this with above (1) is proved. 
\par
Ad(2): 
Since by (1) $d_\lambda= d_{\alpha_\lambda}
$ it follows that $F_\lambda$ defined at the  the proof of (1)
above is the minimal 
conditional expectation.
\par 

For any $a\in {\Cal A}(I),$  we first note that
$$
E(\alpha_{\lambda\bar\lambda}(a))
= E(u^*au)
=u^* E(a) u
=\alpha_{\lambda\bar\lambda}(E(a))
$$
where $u$ is the unitary intertwinner  
transporting $\lambda\bar\lambda$
to $\widehat{\lambda\bar\lambda}$ which is localized on $I_+$. 
we have:
$$
\align
E(F_\lambda(a)) &= 
E(\lambda(R_{\bar\lambda\lambda}) \alpha_\lambda\alpha_{\bar\lambda}(a)
\lambda(R_{\bar\lambda\lambda}^*)) \\
&= \lambda(R_{\bar\lambda\lambda}) E(\alpha_\lambda\alpha_{\bar\lambda}(a))
\lambda(R_{\bar\lambda\lambda}^*)) \\
&=F_\lambda(E(a)) 
\endalign
$$

\enddemo

\qed
\par 
\proclaim{Proposition 3.5}
Suppose $\mu,\lambda$ are representations of ${\Cal B}$ 
localized on $I$. Then:\par
(1) Let $x$ be any subsector of $\alpha_{\mu},$ then
$$
[\alpha_\lambda][x]= [x][\alpha_\lambda];
$$ 
(2)
Let $z$ and $y$ be subsectors of $\alpha_{\lambda}$ and $\tilde\alpha_\mu$
respectively, then
$[z][y]=[y][z].$
\endproclaim
\demo{Proof:}
Ad (1): Let $\epsilon(\lambda,\mu)$ be the braiding operator. 
Then 
$$
\align
\alpha_\mu \alpha_\lambda &= \alpha_{Ad_{\epsilon(\lambda,\mu)} \lambda\mu}
\\
&=\alpha_{Ad_{\epsilon(\lambda,\mu)}} \alpha_\lambda\alpha_\mu
\\
&=Ad_{\epsilon(\lambda,\mu)}\alpha_\lambda\alpha_\mu
\endalign
$$
Now let $v_x\in Hom (x,\alpha_\mu)$ be the isometry such that
$x(\cdot)= v_x^* \alpha_\mu(\cdot) v_x.$ As in the proof of Th. 3.6
on Page 377 of [Xb], it is sufficient to show that
$$
\alpha_\lambda(v_xv_x^*)= \epsilon(\lambda,\mu) v_xv_x^*
\epsilon(\lambda,\mu)^*
$$
Applying $\gamma$ to the above equality, it is   
sufficient to show that
$$
\gamma\alpha_\lambda(v_xv_x^*)= \gamma(\epsilon(\lambda,\mu) v_xv_x^*
\epsilon(\lambda,\mu)^*)
$$
This follows from YBE and BFE as on Page 377 of [Xb].\par
Ad (2): We first prove that
$$
Ad_{\epsilon(\lambda,\mu)} \alpha_\lambda\tilde\alpha_\mu
=\tilde\alpha_\mu \alpha_\lambda
$$ 
Let $u_{\lambda+}$ and $ u_{\mu-}$ be the intertwinners as in the 
definition of $\alpha_\lambda$ and $\tilde\alpha_\mu$  in 
Prop. 3.1. Then the above equality
is equivalent to
$$
\epsilon(\lambda,\mu) u_{\lambda+}^* u_{\mu-}^* =  u_{\mu-}^*
 u_{\lambda+}^* 
$$
Since $\epsilon(\lambda,\mu)= \mu( u_{\lambda+}^*)u_{\lambda+},$
we need to show that
$$
u_{\mu-} \mu( u_{\lambda+}^*)  u_{\mu-}^* =  u_{\lambda+}^*
$$
which follows from the fact that $ u_{\lambda+}^*\in {\Cal B}(I\cup I_+)$
and $\hat \mu$ is localized on $I_-$. \par
The rest of the proof follows by YBE and BFE as on Page 385 of [Xb].
\enddemo
\qed
\par
We note that many properties of relative braidings implicitly used in [Xb] and
further studied in [BE3] can also be proved in our current setting, but
we will not use them in this paper. 
\subheading {\S3.2 Strongly additive pairs}
Let ${\Cal A}$  be a conformal net and 
let ${\Cal B}\subset {\Cal A}$ be a conformal subnet as defined in
\S2.1. Motivated by the proof of Theorem 2.4, we give the following
definition:\par
\proclaim{Definition 3.2}
The pair  ${\Cal B}\subset {\Cal A}$ is said to be strongly 
additive if
$$
{\Cal B}(I_1)\vee {\Cal A}(I_2) = {\Cal A}(I) $$
for any intervals $I,I_1,I_2$ such that 
 $I_1,I_2$ are the connected components
of a set obtained from $I$ by removing an interior point of $I$. 
\endproclaim
Note that the above definition can be generalized to nets of algebras
without conformal invariance, but conformal nets give most interesting
examples of strongly additive pairs, and so we shall consider only conformal
nets in this paper.
Also note that by conformal invariance, it is sufficient to check the
condition in the above definition for a particular $I, I_1,I_2$.
\proclaim{Lemma 3.6}
(1) If the pair  ${\Cal B}\subset {\Cal A}$ is strongly additive, 
then ${\Cal B}$ and ${\Cal A}$ are strongly additive, and
$$ {\Cal B}(I_1)'\cap  {\Cal A}(I_1) = {\Bbb C}, \forall I_1$$
\par
(2) If $G$ is a compact group acting properly on $ {\Cal A}$  and
${\Cal B}$ is the fixed point subnet under the action of $G$, then
the pair  ${\Cal B}\subset {\Cal A}$ is a strongly additive pair 
if and only if either ${\Cal A}$ is strongly additive or
${\Cal B}$ is strongly additive.\par
(3) Let ${\Cal B}\subset {\Cal C}\subset {\Cal A}$ be conformal subnets. 
Then the pair ${\Cal B}\subset {\Cal A}$ is strongly additive iff 
the pairs  ${\Cal B}\subset {\Cal C}$ and ${\Cal C}\subset {\Cal A}$
are strongly additive. 
\endproclaim
\demo{Proof:}
Ad(1):
The first follows trivally from the definition and applying conditional
expectation from ${\Cal A}$ to ${\Cal B}.$ For the second part, choose
$I_2$ which share only one boundary point with $I_1$  and let 
$I$ be the smallest interval containing $I_1\cup I_2$. Then by
the assumption we have
$$
{\Cal B}(I_1) \vee   {\Cal A}(I')=  {\Cal A}(I_2')
$$
Taking the commutants and applying the Haag Duailty, we get
$$
{\Cal B}(I_1)' \cap   {\Cal A}(I)=  {\Cal A}(I_2), 
$$
and so 
$$
{\Cal B}(I_1)' \cap   {\Cal A}(I_1) \subset  {\Cal A}(I_2)
$$ for all $I_2$  sharing only one boundary point with $I_1$.
Choose a sequence of $I_2^{(n)}$ so that $
\cap_n I_2^{(n)} $ is one point, by (c) of Prop. 2.1  we get
$$
{\Cal B}(I_1)' \cap   {\Cal A}(I_1) \subset \cap_n  {\Cal A}(I_2^{(n)})
= {\Bbb C}
$$
\par
Ad(2):
By Theorem 2.4 and (1), it is sufficient to show that if ${\Cal B}$
is strongly additive, then the pair ${\Cal B}\subset {\Cal A}$ 
is strongly additive. 
Let $I$ be an arbitary interval, and $I_1,I_2$ are the connected components
of a set obtained from $I$ by removing an interior point of $I$. 
Then we have the following inclusions:
$$
{\Cal B}(I)\subset {\Cal B}(I_1)\vee {\Cal A}(I_2) \subset {\Cal A}(I) 
$$
where the first inclusion follows from the strong additivity of 
${\Cal B}$. By [I], there exists a closed subgroup $G_1$ of $G$ such that
${\Cal B}(I_1)\vee {\Cal A}(I_2)$ is the fixed point subalgebra of
${\Cal A}(I)$ under the action of $G_1$. It follows that there is 
normal faithful conditional expectation form ${\Cal A}(I)$ to 
${\Cal B}(I_1)\vee {\Cal A}(I_2)$ preserving the vector state
$(\cdot\Omega,\Omega)$. Since $ {\Cal A}(I_2)\Omega$ is dense in 
$H$ by Reeh-Schlieder theorem in Prop. 2.1 it follows that
$${\Cal B}(I_1)\vee {\Cal A}(I_2) = {\Cal A}(I) 
$$ 
by Takesaki's theorem (cf. \S9 of [SS]). 
\par
(3) follows directly from the definitions and applying suitable conditional 
expectations.  

\enddemo
\qed \par
Note that orbifolds as discussed in \S2 are examples of strongly
additive pairs. The following proposition will be used in \S4
to give more strongly additive pairs (cf. the list in \S4.1).
\proclaim{Proposition 3.7}
Let   ${\Cal A}$ be a conformal net and 
${\Cal B}\subset {\Cal A}$ is a conformal subnet. Assume that
${\Cal B}$ is strongly additive and $U(g)\in \vee_I{\Cal B}(I),
\forall g\in {\bold G}$ where $U$ is the representation of 
conformal group ${\bold G}$ as defned in ${\Cal B}$ of \S2.1.
Then the pair ${\Cal B}\subset {\Cal A}$ is strongly additive.
\endproclaim
\demo{Proof:}
Let $I$ be an interval and $I_1,I_2$ the intervals obtained by
removing an interior point of $I$. 
Let $E$ be the unique conditional expectation from $  {\Cal A}(I)$
to $ {\Cal B}(I)$ such that 
$\psi (E(\cdot))= \psi(\cdot)$
where $\psi(\cdot)=(\cdot\Omega,\Omega)$ is the 
normal faithful state on $  {\Cal A}(I)$ and $\Omega$ is the
vacuum vector. Denote by $\Delta_\psi$ and $\Delta$ the modular
operator of  $  {\Cal A}(I)$ and  $  {\Cal B}(I)\vee {\Cal A}(I_1)$  
with respect to $\Omega$. Notice that $Ad_{\Delta_\psi^{it}}$
and  $Ad_{\Delta^{it}}, t\in {\Bbb R}$ induce the same automorphism on 
$ {\Cal B}(I)$ and $ {\Cal B}(I')$, and so
$\Delta_\psi^{it} \Delta^{-it} \in  {\Cal B}(I)'\cap {\Cal B}(I')'$.
By the geometric nature of $\Delta_\psi^{it}$ (cf. Prop. 2.1) and 
our assumption,
$\Delta_\psi^{it} \in  \vee_I{\Cal B}(I)$, and by strong additivity
of $ {\Cal B}$ we have
$\Delta_\psi^{it} \in {\Cal B}(I)\vee {\Cal B}(I'), \forall t\in {\Bbb R}
.$
So $\Delta_\psi^{it} \Delta^{-it}$ commute with $\Delta_\psi^{it'}$ for
all $t,t'\in {\Bbb R}$, hence  
$\Delta_\psi^{it'} $ commutes with  $\Delta^{it}$. It follows that
for any $t'$, 
$Ad_{\Delta^{it}}$ is a one parameter  automorphism of 
$Ad_{\Delta_\psi^{it'}}(
 {\Cal B}(I)\vee {\Cal A}(I_1))$ preserving the vector state
$(\cdot\Omega,\Omega)$.  
By KMS condition (cf. Page 28 pf [SS]) $Ad_{\Delta^{it}}$ is the 
modular automorphism of $Ad_{\Delta_\psi^{it'}}(
 {\Cal B}(I)\vee {\Cal A}(I_1))$ with respect to $\Omega$. Notice that
$\overline{{\Cal A}(I_1)\Omega} =\overline{{\Cal A}(I)\Omega}$
for any interval $I$ by Reeh-Schilider's Theorem, it follows by 
a Theorem of Takesaki (cf. \S9 of [SS]) that
$$
Ad_{\Delta_\psi^{it'}}(
 {\Cal B}(I)\vee {\Cal A}(I_1))=  {\Cal B}(I)\vee {\Cal A}(I_1)
,\forall t'\in {\Bbb R}
$$
Since $\cup_{t'} \Lambda_I(2\pi t')I_1= I$, it follows that
$$
{\Cal B}(I)\vee {\Cal A}(I_1)=  {\Cal A}(I)
$$
which proves the proposition since ${\Cal B}$ is strongly additive. 
\enddemo
\qed
\par
The following  is a generaliztion of Th. 3.3 of [Xb] and is the key result of
\S3:
\proclaim{Theorem 3.8}
Let ${\Cal B}\subset {\Cal A}$ be a strongly addive pair of conformal nets, 
and suppose $\mu,\lambda$ are representations of ${\Cal B}$ 
localized on $I$. \par
(1) 
If $x\in {\Cal A}(I)$ satisfies
$x\lambda (b)= \mu(b)x, \forall b\in {\Cal B}(I),$ 
then $x \alpha_\lambda(a) =  \alpha_\mu(a)x, \forall a\in  {\Cal A}(I).$ 
\par
(2) If $\mu,\lambda$ have finite index, then
$$
\langle \alpha_\mu, \alpha_\lambda\rangle = 
\langle  \mu\bar\lambda, \gamma\rangle
$$
where $\gamma$ is a representation of ${\Cal B}$ unitarily equivalent
to the defining representation of ${\Cal B}$ on the vacuum Hilbert 
space of ${\Cal A}$. 
\par
\endproclaim
\demo{Proof:}
Ad (1): Let $u_1$ (resp. $u_2$) be unitary intertwinner 
in ${\Cal B}(J), J\supset I$
transporting 
$\lambda$ (resp. $\mu$) to $\hat\lambda$ (resp.  $\hat\mu$)
localized on $I_+$ as in the definition of braiding operator (cf. \S2.2). 
Since 
 $x \lambda(a) =  \mu(a)x, \forall a\in  {\Cal B}(I)$
and $x\in {\Cal A}(I),$ it follows that by formula (3) of Prop. 3.1
$u_2 xu_1^* \in  {\Cal A}(J) \cap  {\Cal B}(I)'.$ 
Let us choose $J$ so that $I, I_+$ are the intervals obtained by removing
an interior point of $J$. By the strong additive pair assumption, we have
$$
{\Cal A}(J') \cup  {\Cal B}(I)=  {\Cal A}(I_+') 
$$
and so 
$$
{\Cal A}(J) \cap  {\Cal B}(I)'=  {\Cal A}(I_+) 
$$
It follows that
$$
u_2 xu_1^* \in  {\Cal A}(I_+)
$$
and (1) follows from formula (3) in Prop. 3.1 and 
and locality of the net $ {\Cal A}$. \par
Ad (2):
Since $\lambda$ has finite index, there are two isometries 
$$R_{\lambda \bar\lambda}\in (id, \lambda \bar\lambda), 
R_{\bar\lambda \lambda}\in (id, \bar\lambda\lambda)
$$ in ${\Cal B}(I)$
such that
$$
R_{\lambda \bar\lambda}^*\lambda(R_{\bar\lambda \lambda}) 
=R_{ \bar\lambda\lambda}^*\bar\lambda(R_{\lambda \bar\lambda})
=\frac{1}{d(\lambda)}
$$
Let $x\in Hom (\alpha_\lambda, \alpha_\mu).$ Then
$$
xR_{\lambda \bar\lambda} b =  \mu\bar\lambda(b) xR_{\lambda \bar\lambda}
, \forall b\in {\Cal B}(I)$$
Note that $ xR_{\lambda \bar\lambda}\in  {\Cal A}(I).$
Since $  {\Cal B}(I)\subset  {\Cal A}(I)$ is irreducible by Lemma 3.6,
the vector space $H_{\mu\bar\lambda}:= \{ y\in  {\Cal A}(I)|
yb =\mu\bar\lambda(b) y, \forall b\in  {\Cal B}(I) \}$ 
is a finite dimensional vector space  
with dimension 
$$
\langle \mu\bar\lambda, \gamma \rangle 
$$ by Th. 3.3 (i) of [I]. Note that  the map
$$
x\in Hom(\alpha_\lambda, \alpha_\mu)_{\Cal A}\rightarrow 
xR_{\lambda \bar\lambda}\in H_{\mu\bar\lambda}
$$
is one-to-one by the relations satisfied by $R_{\lambda \bar\lambda},
R_{\bar\lambda \lambda}$. To prove (2) it is sufficient to show that
the following one-to-one map 
$$
y\rightarrow \mu(R_{\bar\lambda \lambda}^*)y
$$
where $y\in  H_{\mu\bar\lambda}$ is a map from
$ H_{\mu\bar\lambda}$ to $ Hom(\alpha_\lambda, \alpha_\mu)_{\Cal A}$.
Note that 
$$ \mu(R_{\bar\lambda \lambda}^*)y 
\lambda(b)= \mu(b)\mu(R_{\bar\lambda \lambda}^*)y $$ 
and $\mu(R_{\bar\lambda \lambda}^*)y\in {\Cal A}(I)$. By (1) 
$$
\mu(R_{\bar\lambda \lambda}^*)y\in Hom(\alpha_\lambda, \alpha_\mu)_{\Cal A}
$$

\enddemo
\qed
\par 
\proclaim{ Corollary 3.9}
With the notations as in Theorem 3.8, and assume that  
$\lambda$ has finite index. 
Then  
(1) $
[\alpha_{\bar\lambda}] = [\bar\alpha_{\lambda}]
$ ;\par
(2)
Let $H_\lambda:=\{x\in {\Cal A}(I)| xb= \lambda(b)x, \forall b\in {\Cal B}(I)
\}$. $H_\lambda$ is called {\it the space of charged intertwinners} 
associated with $\lambda$ as in [LR]. 
Then $H_\lambda= Hom(id, \alpha_\lambda)$ 
and 
$dim H_\lambda= \langle \gamma,\lambda\rangle\leq d_\lambda$. \par
(3) Let ${\Cal A}_f(I)\subset {\Cal A}$ be the subalgebra generated by
${\Cal B}(I)$ and $H_\lambda, \lambda\in S$, where $S$ is a set of 
(DHR) irreducible representations of ${\Cal B}$ with finite 
statistical dimensions, and  is closed under
fusion and conjugation. Then
${\Cal A}_f(I)$ is invariant as a set under the modular automorphism
$AdU(\Lambda_I(t))$ (cf. Prop. 2.1), and 
there exists a unique conformal subnet 
${\Cal L}\subset {\Cal A}$ such that 
$$
{\Cal B}(J)\subset {\Cal L}_f(J)\subset {\Cal A}(J), \forall J\in {\Cal I}
$$
and ${\Cal A}_f(I)= {\Cal L}(I)$. Moreover, the vacuum representation
$H_{\Cal L}$ of  the conformal net ${\Cal L}$ as a representation of 
${\Cal B}$ decomposes as $H_{\Cal L}\simeq\oplus_{\lambda\in S} 
dim{H_\lambda}\lambda$.
\endproclaim
\demo{Proof:}
Ad (1):
The proof is the same as that of [BE1] by Lemma 3.2 and 
(2) of Theorem 3.8 as follows:
$$
\align
\langle \alpha_\lambda, \alpha_\lambda\rangle
&= \langle \lambda\bar\lambda, \gamma \rangle \\
&=\langle  \alpha_{\lambda\bar\lambda}, id \rangle \\
&= \langle  \alpha_\lambda\alpha_{\bar\lambda}, id \rangle \\
&=  \langle  \alpha_\lambda, \bar\alpha_{\bar\lambda}\rangle 
\endalign
$$
Replace $\lambda$ by $\bar\lambda$ we get
$$
\langle \alpha_{\bar\lambda}, \alpha_{\bar\lambda}\rangle= 
\langle  \alpha_{\bar\lambda}, \bar\alpha_{\lambda}\rangle
$$ 
Thus we have
$$
\langle \alpha_{\bar\lambda}, \alpha_{\bar\lambda}\rangle= 
\langle  \alpha_{\bar\lambda}, \bar\alpha_{\lambda}\rangle
=\langle \bar\alpha_\lambda, \bar\alpha_\lambda\rangle 
$$ 
Since $d_{\alpha_\lambda}<\infty, d_{\alpha_{\bar\lambda}} <\infty,$
by Prop. 3.4, 
the above identities imply (1). \par
(2) follows directly from (2) of Theorem 3.8 and (1) of Prop. 3.4.\par
Ad (3):
Let $F_\lambda$ be the minimal conditional expectation as in (3) of Prop. 3.4.
By [L1] we can 
choose a set of isometries $v_{i}\in  Hom(id, \alpha_\lambda), 
1\leq i\leq dimH_\lambda$ such that
$F_\lambda({v_iv_j^*})= \frac{1}{d_\lambda}\delta_{ij}  $.
Note that $E(v_iv_j^*) \in {\Cal B}(I)'\cap {\Cal A}(I)={\Bbb C}$,
and by (3) of Prop. 3.4 we have
$$
E(v_iv_j^*)= F_\lambda (E(v_iv_j^*))
= E( F_\lambda (v_iv_j^*))=\frac{1}{d_\lambda}\delta_{ij}
$$
It follows that the operator $a_\lambda$ as defined on Page 39 of [I]
is the identity operator (our $\lambda$ corresponds to $\xi$ in [I]).
Then the argument on Page 41 of [I] shows that 
$\sigma_t^{\psi E}(v_i)=v_i$ where 
$\psi$ is a dominant weight on ${\Cal B}(I)$, and 
$\sigma_t^{\psi E}$ is the modular 
automorphisms associated with the weight $\psi E$. By Haagerup's Theorem
(cf. Page 156 of [SS]), 
$$
AdU(\Lambda_I(t))(\cdot)=Ad_{u_t} \sigma_t^{\psi E}(\cdot)
$$
where $u_t\in {\Cal B}(I)$. It follows that
${\Cal A}_f(I)$ is invariant as a set under the modular automorphism
$AdU(\Lambda_I(t))$. 
Since $H_\lambda$ is finite dimensional, the rest of the proof is the same as  the proof on Page 18 of [Ls] as follows. Note that $\Lambda_I({\Bbb R})$
is exactly the subgroup of $PSL(2,{\Bbb R})$ which leaves $I$ globally
fixed.  For each $J\in {\Cal I}$, set 
${\Cal L}(J)= Ad_{U(g)}({\Cal A}_f(I))$, where $g\in PSL(2,{\Bbb R}),
gI=J$. It is easy to check that ${\Cal L}(J)$ is independent of the 
choice of $g$ as long as $gI=J$. Note that ${\Cal L}(J)$ verifies locality
since ${\Cal L}(J)\subset {\Cal A}(J)$. To show that  ${\Cal L}$ is a 
conformal net we just have to check the isotony property, namely 
${\Cal L}(J_1)\subset {\Cal L}(J_2)$ if $J_1\subset J_2$. By 
conformal invariance we may assume that
$J_1=I$ and that $J_2=gI$ for some  $g\in PSL(2,{\Bbb R})$, and it is 
sufficient to show that
$Ad_{U(g)}({\Cal A}_f(I))\supset {\Cal A}_f(I)$.
Since $H_\lambda$ is finite dimensional, By the second part of Cor. 19 in [Ls]
$Ad_{U(g)} H_\lambda= z_\lambda(g)^*  H_\lambda$, where
$z_\lambda(g)\in   {\Cal B}(J_2)$ is a unitary defined by 
formula (14) of [Ls]. Hence
$$
\align
Ad_{U(g)}({\Cal A}_f(I)) &= \{  Ad_{U(g)}({\Cal B}(I)), 
 Ad_{U(g)} H_\lambda,\lambda\in S \}''\\
&= 
\{ {\Cal B}(J_2), 
z_\lambda(g)^* H_\lambda,\lambda\in S \}''\\
&= \{ {\Cal B}(J_2), 
H_\lambda,\lambda\in S \}''
\supset \{ {\Cal B}(I), 
H_\lambda,\lambda\in S \}''\\
&={\Cal A}_f(I)  
\endalign
$$
Now let $\Omega$ be the vacuum vector for ${\Cal A}$, since
${\Cal L}$ is a conformal net, the vacuum representation space of
${\Cal L}$ can be identified as $\overline{{\Cal A}_f(I)\Omega}$
by Reeh-Schilider theorem in Prop. 2.1. For each $\lambda\in S$,
we choose isometries $v_{\lambda,i}, 1\leq i\leq dimH_\lambda$
as in the beginning of the proof of (3) (we add a 
subscript $\lambda$ to emphasize its dependence on $\lambda$).
Then the set consisting of 
$\sum_{\lambda\in S} v_{\lambda,i}^* x_{\lambda,i}, x_{\lambda,i}\in 
{\Cal B}(I)$ where the sum is a finite sum is a dense subalgebra of
${\Cal A}_f(I)$. Note that the 
space $X_{\lambda,i}:=v_{\lambda,i}^*\overline{{\Cal B}(I)\Omega}$
is invariant under the action of ${\Cal B}$, and  the restriction of
 ${\Cal B}$ to this space is a representation of  ${\Cal B}$
unitarily equivalent to $\lambda$. Also note that
$X_{\lambda,i}\bot  X_{\lambda',i'}$ if $(\lambda,i)\neq (\lambda',i') 
$ since $E( v_{\lambda,i}  v_{\lambda',i'}^*)= \frac{1}{d_\lambda}
\delta_{\lambda\lambda'}
\delta_{ii'}$. Hence  $\overline{{\Cal A}_f(I)\Omega}$ is a direct sum of
$X_{\lambda,i}$, and  this proves the last part of (3).
\enddemo
\qed
\par
We will call the conformal net ${\Cal L}$ constructed in (3) of Cor. 3.9
{\it the conformal subnet of ${\Cal A}$ generated by ${\Cal B}$ and 
charged intertwinners associated with the set $S$}.
\heading \S4. Applications \endheading 
\subheading{\S4.1 Conformal nets with central charge 1}
The irreducible representations of Virasoro algebra with central 
charge 1 are classified 
as follows. For each $n\geq 0$, there is an irreducible representation
with lowest weight $n$, and such representation is denoted by $L(1,n)$
as in [D]. Here $1$ is the central charge. When $n=m^2, 2m\in {\Bbb Z}$,
$L(1,n)$ will be called the {\it degenerate representations} due to
the degenerate nature of certain Verma modules. The vacuum representation
is $L(1,0)$. All $L(1,n)$ can be ``exponentiated'' to give irreducible
projective representations of $Diff(S^1)$ (cf. [GW]),
where $Diff(S^1)$ is the group of smooth differmorphisms of $S^1$. 
 On the vacuum representation
$L(1,0)$,
one can define a conformal net ${\Cal A}_{c=1}$ as in \S3 of [FG], 
called Virasoro net with central charge $1$. \par
Let ${\Cal A_{SU(2)_1}}$ be the conformal net associated with 
loop group $LSU(2)$ at level $1$. The adjoint action of 
group $SO(3)$ acts properly on ${\Cal A_{SU(2)_1}}$ in the
sense of \S2.3, and  the fixed point is identified  as ${\Cal A}_{c=1}$, the 
Virasoro net with central charge $c=1$, in [R3].\par
It follows from Theorem 2.4 that ${\Cal A}_{c=1}$ is strongly additive. 
It is already pointed out in [R3] that the statistical dimension 
$d_{L(1,m^2)}$ of
$L(1,m^2)$ is $2m+1$ when $m$ is a non-negative integer, and the
fusion rules among $L(1,m^2)$ are the same as representation 
rings of $SO(3)$.  The following lemma
generalize this to the case when $2m$ is a non-negative integer:
\proclaim{4.1 Lemma}
Assume that $2m$ are non-negative integers.
Then $d_{L(1,m^2)}= 2m+1$ and the fusion ring generated by
$L(1,m^2)$ is  isomorphic  to the representation ring of $SU(2)$.
\endproclaim
\demo{Proof:}
Consider the following conformal inclusion (cf. \S3.1 of [Xj]) 
$$
LSU(2)_1\times LU(1)_2 \subset LU(2)_1
$$
The group $SU(2)$ acts properly on the  
net ${\Cal A}_{U(2)_1}$ with fixed point
net ${\Cal A}_{c=1}\times {\Cal A}_{U(1)_2}$. We note that 
the  
net ${\Cal A}_{U(2)_1}$ is not local, but satisfies the twisted duality
(cf. \S15 of [Wa]). So  Th. 3.6 in [DR2] apply in this case, and the fusion
of those    
irreducible representations of ${\Cal A}_{c=1}\times {\Cal A}_{U(1)_2}$ 
appearing in  ${\Cal A}_{U(2)_1}$ are given by the fusion ring of finite
dimensional representations of $SU(2)$. Since  irreducible covariant 
representations of  
${\Cal A}_{U(1)_2}$ have statistical
dimension equal to $1$ and generate an abelian group ${\Bbb Z}_2$
under the fusion,  the lemma follows. 
\enddemo
\qed
\par
Let ${\Cal A}$ be a conformal net. Following [KL], we say that ${\Cal A}$
is a {\it diffeomorphism covariant net}  if there exists
a unitary projective representation $U$ of $Diff (S^1)$ on $H$ extending the 
unitary representation of $PSL(2,{\Bbb R})$ such that 
$$
U(g){\Cal A}(I)U(g^*)= {\Cal A}(g.I), g\in Diff (S^1), I\in {\Cal I}
$$  
We say that ${\Cal A}$ is a {\it conformal net with central charge 1
} if  ${\Cal A}$ is a  diffeomorphism covariant net containing  
${\Cal A}_{c=1}$ as a conformal subnet such that  
$U(Diff(I))'' =({\Cal A}_{c=1})(I), 
\forall I\in {\Cal I}$,
where $Diff(I)$ denotes the group of smooth diffemorphisms $g$ of $S^1$
satisfying  $g(t)=t, t\in I'$.\par
Let us describe the  known list of
such nets. Let $G$ be a closed group of $SO(3)$. Such groups are
well known to be of $A-D-E$ groups corresponding to affine $A-D-E$
graphs. 
Let $\hat G$ be two-fold covering group of $G$ in $SU(2)$. 
We note the Perron-Frobenius eigenvectors given on Page 14
of [GHJ] are the dimensions of the irreducible representations of
the two-fold covering group $\hat G$. 
Since ${\Cal A}_{c=1}$ can be identified with 
${\Cal A}_{SU(2)_1}^{SO(3)}$, ${\Cal A}_{SU(2)_1}^{G}$ is a 
conformal net with central charge $c=1$. 
The remaining two cases are 
${\Cal A}_{U(1)_{2n}}$ and its ${\Bbb Z}_2$ orbifold  
${\Cal A}_{U(1)_{2n}}^{{\Bbb Z}_2}$ as studied in [Xo], where
$n$ is not the square of an integer. So the known list of
conformal nets with central charge 1 is:
$$
{\Cal A}_{SU(2)_1}^{G}, {\Cal A}_{U(1)_{2n}}, 
{\Cal A}_{U(1)_{2n}}^{{\Bbb Z}_2} \tag * 
$$
where $G$ is a closed subgroup of $SO(3)$ and $n$ is not 
the square of an integer. \par
It has been conjectured (cf [DVVV]) that the list (*) exhausts all
conformal theories with central charge $1$. \par
When $G$ is a finite group, 
${\Cal A}_{SU(2)_1}^{G}$ is absolutely rational by Prop. 2.2. \par
$
{\Cal A}_{U(1)_{2n}}, 
{\Cal A}_{U(1)_{2n}}^{{\Bbb Z}_2} 
$ are also  absolutely rational and all irreducible representations
are obtained in [Xo]. The  irreducible representations of
$
{\Cal A}_{U(1)_{2n}}$ will be denoted by $\pi_i, i\in {\Bbb Z}_{2n}$. They
generate a fusion ring isomorphic to ${\Bbb Z}_{2n}$. \par
When $G=U(1)$, ${\Cal A}_{SU(2)_1}^{G}$ is the net corresponding
to the  Heisenberg group, denoted by $H(1)$. 
$H(1)$ is the set $C^\infty(S^1,{\Bbb R})\times S^1$ with multiplication  
defined  by (cf. \S9.5 of [PS])
$$(f_1, x_1) 
\cdot( f_2, x_2)= (f_1+f_2,  e^{\int_{S^1}f_1f_2'}x_1x_2)
$$ 
Note that $(x,1)$ where $x\in {\Bbb R}$ is considered as a constant map
is in the center of  $H(1)$. For each real number $q$, there is an irreducible representation 
of  $H(1)$ denoted by $F_q$ where $(x,1)$ acts on  $F_q$ as
$(x,1)\rightarrow qx$, and these are all the irreducible representations
of $H(1)$ (cf. Prop. 9.5.10 of [PS]). The net  ${\Cal A}_{SU(2)_1}^{U(1)}$
is related to  $H(1)$ as follows.  $F_0$ is the vacuum representation
of ${\Cal A}_{SU(2)_1}^{U(1)}$, and 
${\Cal A}_{SU(2)_1}^{U(1)}(I)= \pi_{F_0}(C_0^\infty(I,{\Bbb R}))''$ 
where $C_0^\infty(I,{\Bbb R})$ is the set of smooth maps from
$I$ to ${\Bbb R}$ which vanishes on the boundary, and is considered
as a subspace of $C^\infty(S^1,{\Bbb R})$.
Each  $F_q$ is also  an irreducible representation 
of  ${\Cal A}_{SU(2)_1}^{U(1)}$. 
The net  ${\Cal A}_{SU(2)_1}^{U(1)}$ was  studied in [BMT]. 
\par

Some decompositions of the vacuum representation of the nets in the 
list above when restricting to  ${\Cal A}_{c=1}$ 
are also known (cf. Prop. 2.2, Th. 2.7 and Th. 2.9 of [D] where 
our $2n$ corresponds to $n$ in [D]) as follows:\par 
If $n$ is not the square of an integer, then
$$
\align
H_{{\Cal A}_{U(1)_{2n}}}&= (\oplus_{p\geq 0}L(1,p^2))\oplus(\oplus_{m>0}
2L(1,m^2n))  \\
H_{{\Cal A}_{U(1)_{2n}}^{{\Bbb Z}_2} }&= 
(\oplus_{p\geq 0}L(1,4p^2))\oplus(\oplus_{m>0}
L(1,m^2n)) 
\endalign
$$
If $n=k^2$ where $k$ is a nonnegative integer then
$$
H_{{\Cal A}_{U(1)_{2n}}}= \oplus_{m\geq 0}\oplus_{0\leq p\leq {k-1}}
(2m+1)L(1,(mk+p)^2)
$$
When $G=U(1)$ or $G=D_\infty$(infinite dihedral group), we have:
$$
\align
H_{{\Cal A}_{SU(2)_1}^{U(1)}}&= \oplus_{p\geq 0}L(1,p^2)\\
H_{{\Cal A}_{SU(2)_1}^{D_\infty}}&= \oplus_{p\geq 0}L(1,4p^2)
\endalign
$$
Recall that $F_q$ is the irreducible representation 
of  ${\Cal A}_{SU(2)_1}^{U(1)}$  corresponding to 
an irreducible representation of  $H(1)$ labeled by a real number $q$.
The decompositions of $F_q$ with respect to ${\Cal A}_{c=1}$ is also
well known to be (cf. [D]) the following: \par
If $q=\frac{p^2}{4}$ for some non-negative integer $p$, then
$$
F_q= \oplus_{-\frac{1}{2}p\leq m\leq \frac{1}{2}p, m+\frac{1}{2}p
\in {\Bbb Z}}{L(1,m^2)}
$$
If $4q$ is not the square of an integer, then
$F_q= L(1,\sqrt q)$. 
\par

We note that by definition if   ${\Cal A}$ is a conformal
net  with central charge 1, then the pair ${\Cal A}_{c=1}\subset 
{\Cal A}$ satisfies the condition of Prop. 3.7, and hence is
a strongly additive pair in the sense of definition 3.2. 
By Lemma 4.1, the principal graph of $L(1,\frac{1}{4})$ is $A_\infty$, 
and $d_{L(1,\frac{1}{4})}=2$.  It follows from Prop. 3.4 
that $\alpha_{L(1,p^2)}
=d_{L(1,p^2)}= 2p+1$ for all non-negative integer $2p$. Since 
$\alpha_f:=\alpha_{L(1,\frac{1}{4})}$ has minimal index $4$, its principal
graph are determined in [GHJ] and [Po]. In the following
Lemma we list properties of   $\alpha_f$:
\proclaim{Lemma 4.2}
The possible principal graphs of $\alpha_f$ are given by the A-D-E graphs
on Page 19 of [GHJ]. More precisely:\par
(1) $[\alpha_f]=[\bar \alpha_f]$; \par
(2) If $\alpha_f$ is irreducible, then its principal graph is given by 
$D,E$ graph on Page 19 of [GHJ] or $A_\infty$ and $D_\infty$ on 
Page 217 of [GHJ].  \par
(3) If $\langle\alpha_f, \alpha_f\rangle=4$, then its   principal graph
is $A_1^{(1)}$. \par
(4) If $\langle\alpha_f, \alpha_f\rangle=2$, then the principal graph is 
given by either $A_{\infty,\infty}$ or $A_n^{(1)}$.
\endproclaim
\demo{Proof:}
(1) follows from Cor. 3.9. (2), (3) and (4) follows from \S4 of [Po1]. 
%
\enddemo
\qed
\par
As a warm up exercise , let us work out $\alpha_f$ for the list (*). We first prove the following
Lemma which will also be used in \S4.2:
\proclaim{Lemma 4.3}
Every covariant representation of  ${\Cal A}_{SU(2)_1}^{U(1)}$ 
is a direct sum of of irreducible representations, and every irreducible
covariant representation of $ {\Cal A}_{SU(2)_1}^{U(1)}$ 
is isomorphic to some $F_q$.
\endproclaim
\demo{Proof:}
Let $\pi$ be a covariant representation of  ${\Cal A}_{SU(2)_1}^{U(1)}$ . 
Recall the product rule in $H(1)$: for any 
$f_i:S^1\rightarrow {\Bbb R},i=1,2,$
$$
(f_1,x_1)(f_2,x_2)= (f_1+f_2, x_1x_2e^{i\int_{S^1}f_1f_2'})
$$  
We will write $(f_i,1)$ simply as $f_i$ in the following.\par  
Let $I_i\in {\Cal I}, 1\leq i\leq n$ be an open covering of $S^1$ and $\phi_i
,1\leq i\leq n$ a partition of unity such that $supp(\phi_i)\in I_i, 
1\leq i\leq n$. If $f_:S^1\rightarrow {\Bbb R}$,
we assume that in $H(1)$:
$$
f= \prod_k(f\phi_k,1)  C(f,\phi)
$$
where $C(f,\phi)\in {\Bbb C}$ is a phase coming from the product rules.
we have $\pi_{F_0}(f\phi_k)\in {\Cal A}_{SU(2)_1}^{U(1)}(I_k)$, and we define
$$
\pi(f):= \prod_k\pi_{I_k} (\pi_{F_0}(f\phi_k))C(f,\phi)
$$
It is routine to check that $\pi(f)$ is independent of the choices of
open coverings and partition of unity by using the product rules
in  $H(1)$ and isotony, and it gives a  representation of
$H(1)$ with positive energy. The Lemma now follows from Prop. 9.5.10
of [PS] and its proof. 
\enddemo
\qed
\par
\proclaim{Lemma 4.4}
(1) If ${\Cal A}= {\Cal A}_{c=1}$, then the principal graph of
$\alpha_f$ is $A_\infty$;\par
(2) If ${\Cal A}= {\Cal A}_{U(1)_{2k^2}}$, then the principal graph of
$\alpha_f$ is $A_{2k-1}^{(1)}$;\par
(3)  If ${\Cal A}= {\Cal A}_{U(1)_{2n}}$ where $n$ is not the square of an 
integer, 
or  ${\Cal A}= {\Cal A}_{SU(2)_1}^{U(1)}$  then the principal graph of
$\alpha_f$ is $A_{\infty,-\infty}$;\par
(4) If ${\Cal A}= {\Cal A}_{U(1)_{2k^2}}^{\Bbb Z_2}$, 
then the principal graph of
$\alpha_f$ is $D_{k}^{(1)}$;\par
(5) If ${\Cal A}= {\Cal A}_{U(1)_{2n}}$ where $n$ is not the square of an 
integer, 
or  ${\Cal A}= {\Cal A}_{SU(2)_1}^{D_\infty}$ , then the principal graph of
$\alpha_f$ is $ D_{\infty,-\infty}$;\par
(6)  If ${\Cal A}= {\Cal A}_{SU(2)_1}^{E_i}, i=6,7,8$, 
then the principal graph of
$\alpha_f$ is $E_i$.
\endproclaim
\demo{Proof:}
(1) follows from Lemma 4.1. If $ {\Cal A}= {\Cal A}_{U(1)_{2k^2}}$, it follows
from the branching rules, Lemma 4.1 and Prop. 3.1 that
$\alpha_f$ is localized on $I$, and by (2) of Th 3.8 
$\langle \alpha_f, \alpha_f\rangle= 4$ when $k=1$ and
$\langle \alpha_f, \alpha_f\rangle= 2$ when $k>1$. 
So (2) is proved for $k=1$ by 
Lemma 4.2. We note that 
when $k=1$, ${\Cal A}_{U(1)_2}$ can be identified with ${\Cal A}_{SU(2)_1}$, 
and it is easy to check that
$[\alpha_f]= 2[\tau]$, where $\tau$ is the irreducible representation of
${\Cal A}_{SU(2)_1}$ which is not the vacuum representation and
$[\tau]^2=[1]$. \par
When $k>1$, $[\alpha_f]= [\sigma_1]+[\sigma_2]$ where 
$\sigma_1,\sigma_2$ are representations of ${\Cal A}_{U(1)_{2k^2}}$
which are classified, and by (1) of Lemma 3.3 we must have
$[\alpha_f]= [\pi_{k}]+ [\pi_{-k}]$ and (2) follows. If 
 ${\Cal A}= {\Cal A}_{U(1)_{2n}}$ where $n$ is not the square of an 
integer, 
or  ${\Cal A}= {\Cal A}_{SU(2)_1}^{U(1)}$, it follows from  the branching rules, Lemma 4.1 and  (2) of Th 3.8  that 
$\langle \alpha_{L(1,m^2)},  id \rangle=1,$
$\langle \alpha_{L(1,m^2)},  \alpha_{L(1,m^2)} \rangle=2m+1, 
\forall m\in {\Bbb N},$
$\langle \alpha_f, \alpha_f\rangle=2.$ 
So  $[\alpha_f]= [\sigma_1]+[\sigma_2]$ and
$d_{\sigma_1}= d_{\sigma_2}=1$. By lemma 4.2 $[\alpha_f]=[\bar \alpha_f]$, so
either $[\bar \sigma_1]=[\sigma_2]$ or $[\bar \sigma_i]=[\sigma_i], i=1,2$.
If $[\bar \sigma_i]=[\sigma_i], i=1,2$, then
$[\sigma_i^2]=[1], i=1,2$ and 
$$
[\alpha_{L(1,1)}]=  [\alpha_f^2]-[1]
=[\sigma_1\sigma_2]+ [\sigma_1\sigma_2]+[1]
$$ 

Note that by (1) of Prop. 3.5 we have
$$
[\sigma_1\alpha_f]= [\sigma_1\sigma_2]+[1]=[\alpha_f\sigma_1]= 
[\sigma_2\sigma_1]+[1]
$$
Hence $[\sigma_1\sigma_2]=[\sigma_2\sigma_1]$ implying that
$ \langle\alpha_{L(1,1)},  \alpha_{L(1,1)} \rangle= 5$,  which contradicts
$$\langle \alpha_{L(1,1)},  \alpha_{L(1,1)} \rangle=3$$ So we have
$[\alpha_f]= [\sigma_1]+[\sigma_2], [\bar \sigma_1]=[\sigma_2]$, and it follows that
$$[\alpha_{L(1,m^2)}]= \sum_{-m\leq k\leq m}[{\sigma_1^{2k}}]$$ Since
$\langle \alpha_{L(1,m^2)},  1 \rangle=1$, it follows that 
$[\sigma_1^k]\neq 
[1], 
\forall k,$ and the principal graph of $\alpha_f$ is  $A_{\infty,-\infty}$.
\par
When  ${\Cal A}= {\Cal A}_{U(1)_{2k^2}}$ or 
${\Cal A}= {\Cal A}_{SU(2)_1}^{E_i}, i=6,7,8,$
$\alpha_f$ is irreducible and
a  representation of ${\Cal A}$, denoted by $f_1$. We will call this 
representation the { \it vector representation}.   
Since ${\Cal A}$ is completely
rational, it follows that the graph of $\alpha_f$ must be finite and hence 
must be finte $D,E$ type. Consider the following inclusions of
confomal nets:
$${\Cal A}_{c=1}\subset {\Cal A}\subset {\Cal A}_{SU(2)_1}$$
consider the induction of $f_1$ (as a local representation of
$ {\Cal A}$ ) from ${\Cal A}$ to $ {\Cal A}_{SU(2)_1}$, denoted by
$\alpha_{f_1}$ which by 
(3) of Lemma 3.3, is the same as 
$\alpha_{L(1,\frac{1}{4})}^{{\Cal A}_{c=1}\rightarrow  {\Cal A}_{SU(2)_1}}$. 
But from the proof above
$[\alpha_{L(1,\frac{1}{4})}^{{\Cal A}_{c=1}
\rightarrow  {\Cal A}_{SU(2)_1}}]=2[\tau]$, where  
$\tau$ is the irreducible representation of
${\Cal A}_{SU(2)_1}$ which is not the vacuum representation and
$[\tau]^2=[1]$.
For any
irreducible representations $\lambda\prec f_1^{2n}$ for 
some $n\in {\Bbb N}$, let $m_\lambda:=\langle \alpha_\lambda, 1\rangle$, 
it follows that $m_\lambda=d_\lambda$, and by (2) of Th. 3.8, 
$m_\lambda$ is the multplicity of $\lambda$ which appears in the 
vacuum representation of $ {\Cal A}_{SU(2)_1}$, and is equal to the 
dimension of the representation  of the 
corresponding $D-E$ groups since ${\Cal A}$ is the
fixed point net of  ${\Cal A}_{SU(2)_1}$ under the action of such group.
Note that the Perron-Frobenius
eigenvectors listed on Page 14 of [GHJ] are the statistical dimensions
associated with the corresponding representations. Hence the graph
is uniquely determined by the corresponding groups to be those of
(4) and (6).  \par
(5) follows by inspecting the branching rules and using 
Th. 3.8 as we have done in proving (1)-(3).
\enddemo
\qed
\par
As a by-product of the proof above, we have 
the following proposition which contains the main result of [R1]: 
\proclaim{ Proposition 4.5}
Let $2m$ be a non-negative integer, and $n\geq 0$.  Then
$$
[L(1,m^2)][L(1,n)]= \sum_{-m\leq k\leq m, k+m\in {\Bbb Z}} [L(1, (k+\sqrt{n})^2]
$$
\endproclaim
\demo{Proof:}
By Lemma 4.1 it is sufficient to consider the case when
$L(1,n)$ is generic, i.e., when $4n$ is not the square of an integer. 
Condsider the inclusion ${\Cal A}_{c=1}\subset {\Cal A}_{SU(2)_1}^{U(1)}$. 
It is a strongly additive pair. As in the proof of (3) in Lemma 4.4, 
by the branching rules, Lemma 4.1 and Prop. 3.1, 
$\alpha_f$ is a (DHR) representation of ${\Cal A}_{SU(2)_1}^{U(1)}$, and
$[\alpha_f]=[\sigma_1]+[\sigma_2], [\bar\sigma_1]=[\sigma_2]$. 
By Lemma 4.3, $[\sigma_1]=[F_q]$ for some real $q$, and since
by (1) of Lemma 3.3 
$$
\langle \alpha_f, F_q\rangle \leq \langle L(1,\frac{1}{4}), \gamma F_q\rangle
,$$
inspecting the branching rules we must have 
$[\alpha_{L(1,\frac{1}{4}}]= [F_{\frac{1}{2}}]+  [F_{-\frac{1}{2}}]$. 
It follows by Lemma 4.1 
$$
[\alpha_{L(1,m^2)}] = \sum_{-m\leq k\leq m, k+m\in {\Bbb Z}} [F_k]
$$
Note that $[\gamma(F_{\sqrt{n}})] =[L(1,n)]$ by the branching rules and 
Prop 3.1 of [LR]. 
We have:
$$
\align
[L(1,m^2)][L(1,n)] &=  [L(1,m^2)][\gamma F_{\sqrt{n}}] \\
&= [\gamma \alpha_{L(1,m^2)} F_{\sqrt{n}}]\\
&= [\gamma]\sum_{-m\leq k\leq m, k+m\in {\Bbb Z}} [F_kF_{\sqrt{n}}]  
\\
&=  [\gamma]\sum_{-m\leq k\leq m, k+m\in {\Bbb Z}} [F_{k+\sqrt{n}}]
\\
&=  \sum_{-m\leq k\leq m, k+m\in {\Bbb Z}} [L_{(k+\sqrt{n})^2}]
\endalign
$$
\enddemo
\qed
\subheading{\S4.2 Classifications}
As in [KL], two conformal nets ${\Cal A}$ and ${\Cal B}$ are said to be
isomorphic, denoted by
$ {\Cal A}\simeq {\Cal B}$, 
if there is a uniatry operator $U: {\Cal H}_{\Cal A}\rightarrow
{\Cal H}_{\Cal B}$ such that 
$ U^*{\Cal B}(J)U = {\Cal A}(J), \forall J\in {\Cal I}, 
U\Omega_{\Cal A}=\Omega_{\Cal B},
$  where $\Omega_{\Cal A}$ and  $\Omega_{\Cal B}$ are the vacuum vectors of
${\Cal A}$ and ${\Cal B}$ respectively. \par
Notice by the proof of Lemma 4.4, we can infer the type of principal
graphs for $\alpha_f$ from branching rules and $\gamma$ from 
Th. 3.8.  On the other hand, for the same reason 
due to Th. 3.8, we can deduce information about $\gamma$ from the
type of principal graph of $\alpha_f$. This is the basic strategy which 
we will follow to classify conformal nets with central 
charge $c=1$. This will work out under the following spectrum 
condition:
\proclaim{Spectrum condition:}
We say that a conformal net  with central 
charge $c=1$ verfies the spectrum condition if  a degenerate
represenation of the Virasoro net other than the vacuum representation 
must appear in  the vacuum representation of 
${\Cal A}$ if ${\Cal A}\neq {\Cal A}_{c=1}$.
\endproclaim 
\proclaim{Theorem 4.6}
If  a conformal net  ${\Cal A}$  with central 
charge $c=1$  verfies the spectrum condition,  then   
${\Cal A}$ is isomorphic to one of the nets in the list (*).
\endproclaim
The proof is divided into the following steps: 
\subheading{4.2.1 Discrete case: Full spectrum}
In this section we assume that \par
$\gamma=\oplus_{0\leq m} (2m+1) L(1,m^2)$.
We'd like to show that ${\Cal A}\simeq {\Cal A}_{SU(2)_1}$. 
This is essentially an application of resconstruction theorem of 
Doplicher-Roberts (cf. [DR1-2]).We give a sketch of the proof and
refer to [DR1-2] for details. First by our assumption we can
assume that ${\Cal A}$ and ${\Cal A}_{SU(2)_1}$ acting on the same
Hilbert space, and ${\Cal A}_{c=1}$ is a common conformal subnet. 
Let $\Delta:=\{ L(1,p^2), p\in {\Bbb Z} \}$.  We note that
$L(1,1)$ has permutation symmetry and satisfies special 
conjugate property, and $\Delta$ is generated by  $L(1,1)$, and
is specially directed as defined on Page 98 of [DR1]. 
Let us fix an interval $I$. Choose charged intertwinners 
$\psi^{\Cal A}\in {\Cal A}(I) , \psi\in {\Cal A}_{SU(2)_1}(I)$ for
the set $\Delta$ (it is enough to choose  charged intertwinners for
$L(1,1)$). Denote by ${\Cal U}_{\Cal A}$ (resp. ${\Cal U}$)
the $C^*$ algebra generated by ${\Cal U}_ {{\Cal A}_{c=1}}$ and
$\psi^{\Cal A}$ (resp. $\psi$).   
Note that ${\Cal U}_{\Cal A}\cap {\Cal U}_ {{\Cal A}_{c=1}}'=
{\Bbb C}$.  By Page 93-4 of [DR1] there exists an epimorphism
$\phi: {\Cal U}\rightarrow {\Cal U}_{\Cal A}$ such that
$\phi=id$ on ${\Cal U}_ {{\Cal A}_{c=1}}$. Using $\phi$ one can define 
an action of $SO(3)$ 
which  commutes with ${\Cal U}_ {{\Cal A}_{c=1}}$. 
One checks that $\phi$ commutes with the adjoint action of $SO(3)$, and so
ker$\phi$ is  $SO(3)$ invariant. Then the argument of Page 95 of [DR1]
shows that $ker\phi=\{0\}$, and so $\phi$ is an isomorphism. 
Now define a unitary operator $U$ on $H$ by
$Um\Omega= \phi(m) \Omega$ where $\Omega$ is the vacuum vector. Then we have
$$
U{\Cal U}U^* = {\Cal U}_{\Cal A}
$$
and $U$ commutes with  ${\Cal U}_ {{\Cal A}_{c=1}}$. Passing to the
von Neumann algebra generated by ${\Cal U}$ and ${\Cal U}_{\Cal A}$,
we have 
$U{\Cal A}(I)U^*= {\Cal A}_{SU(2)_1}(I), U\Omega=\Omega.$
Since $U$ commutes with  ${\Cal U}_ {{\Cal A}_{c=1}}$, $U$ commutes with
$PSL(2,{\Bbb R})$ by the strong additivity of ${\Cal A}_{c=1}$, and it follows that 
$$
U{\Cal A}(J)U^*= {\Cal A}_{SU(2)_1}(J), \forall J\in {\Cal I}, U\Omega=\Omega.
$$
\subheading{\S4.2.2 General discrete case}
In this section we assume that $$\gamma=\oplus_{p\geq 0} m_pL(1,p^2)$$ 
Notice that by (2) of Cor. 3.9 $m_p\leq {2p+1}$. By Lemma 4.1 and Prop. 3.1,
$\alpha_{L(1,p^2)}$ is localized. 
Consider the following set of  representations 
$\Delta:\{\lambda\prec
 \alpha_{L(1,1)}^{n},n \in {\Bbb N}\}$ of ${\Cal A}$. 
 Since $\epsilon(\alpha_{L(1,1)}, \alpha_{L(1,1)})= \epsilon(L(1,1),L(1,1))$
by (2) of Lemma 3.3, the set $\Delta$ has permutation symmetry, 
and ${\Cal A}$ is strongly additive, 
we can apply
Doplicher-Roberts reconstruction as in Prop. 3.9 of [Mu]  
to obtain conformal net ${\Cal C}$ such that 
${\Cal A}_{c=1}\subset {\Cal A} \subset {\Cal C}$.
Consider the conformal nets
$
{\Cal A}_{c=1} \subset {\Cal C}
$.
We claim that $m_{p}^{\Cal C}= 2p+1$ in this case, where
$m_{p}^{\Cal C}$ is the multiplicity of $L(1,p^2)$ which appears in the
vacuum representation of ${\Cal C}$. 
We have 
$m_{p}^{\Cal C} 
\geq\sum_{\lambda} \langle  \alpha_{L(1,p^2)}, \lambda\rangle$,
where the sum is over those irreducible representations 
$\lambda$ of ${\Cal A}$ which appears in the  
vacuum representation of ${\Cal C}$. We note that 
$\sum_{\lambda} \langle  \alpha_{L(1,p^2)}, \lambda\rangle$ is completely 
determined by the principal graph of $\alpha_f$, which corresponds to
A-D-E groups. So the number 
$\sum_{\lambda} \langle  \alpha_{L(1,p^2)}, \lambda\rangle$ 
depends only the type of A-D-E graph associated with  $\alpha_f$. 
Since all types of such graphs have appeared in Lemma 4.4, and in each case,
it is easy to check the number 
$\sum_{\lambda} \langle  \alpha_{L(1,p^2)}, \lambda\rangle$ is $2p+1$, since
for ${\Cal A}_{c=1}\subset {\Cal A}_{SU(2)_1}^G$, where $G$ is a closed
subgroup of $SO(3)$, the above reconstruction give us ${\Cal A}_{SU(2)_1}$. 
It follows 
that   $m_{p}^{\Cal C}\geq  2p+1$, and by (2) of Cor. 3.9
$m_{p}^{\Cal C}= 2p+1$. 
Now considered the conformal subnet ${\Cal D}\subset {\Cal C}$
generated by ${\Cal A}_{c=1}$ and charged intertwinners for 
$\{ L(1,p^2): \forall p\geq 0 \}$ whose existence is shown by Cor. 3.9.
${\Cal D}$ now has full spectrum as in \S4.2.1, and so
${\Cal A}_{c=1} \subset {\Cal A} \subset
{\Cal D}\simeq {\Cal A}_{SU(2)_1}$. 

Since ${\Cal A}_{c=1}$ is the fixed point net of ${\Cal A}_{SU(2)_1}$
under the action of $SO(3)$, by [I] there exists a closed subgroup
$G_1$ of $SO(3)$ such that
${\Cal A}(I)$ is the fixed point subalgebra of  ${\Cal A}_{SU(2)_1}(I)$
under the action of $G_1$. Since $G_1$ commutes with $PSL(2,{\Bbb R})$,
${\Cal A}$ is the fixed point net of ${\Cal A}_{SU(2)_1}$
under the action of $G_1$. Since $G_1$ is classified as one of the
$A-D-E$ groups, it follows that  ${\Cal A}$ is in the list (*).
\subheading{4.2.3 
${\Cal A}_{c=1}\subset {\Cal A}_{U(1)_{2n}} \subset {\Cal A}$ or 
${\Cal A}_{c=1}\subset {\Cal A}_{SU(2)_1}^{U(1)} \subset {\Cal A}$ case
}
Let us first consider the case when 
${\Cal A}_{c=1}\subset {\Cal A}_{U(1)_{2n}} \subset {\Cal A}$. \par 
Recall that ${\Cal A}_{U(1)_{2n}}$ 
is completely rational, and its representations
are labeled by $\pi_i, 0\leq i\leq 2n-1$, with conformal dimensions
(the eigenvalues of the action of rotations)
$\frac{i^2}{4n^2}+{\Bbb N}$ and statistical
dimension $1$. The fusion ring generated by  $\pi_i$ is ${\Bbb Z}_{2n}$.
Since ${\Cal A}_{U(1)_{2n}}\subset  {\Cal A}$
is a strongly additive pair by (3) of Lemma 3.6, the vacuum representation
$H_{\Cal A}$ decomposes into representation of   ${\Cal A}_{U(1)_{2n}}$ 
as $H_{\Cal A}=\oplus_{i} m_i \pi_i$ with $m_i\leq 1$ by (2) of Cor. 3.9.
Also not that if $m_i=1$, then again by (2) of Cor. 3.9 $
m_i=1=\langle 1, \alpha_i\rangle$, hence  $[\alpha_i]=[1]$ since
$d_{\alpha_i}=d_i =1$. So the set of $H_i$ which appears in the 
decompositions of $H_{\Cal A}$ form an abelian subgroup of ${\Bbb Z}_{2n}$. 
Let $i_0>0$ be the generator of this abelian subgroup. Then there is a
positive integer $k$ such that $ki_0= 2n$. On the other hand, we must have
$\frac{i_0^2}{4n^2}\in {\Bbb Z}$, and we assume that $k_0$ is a positive 
integer such that $i_0^2=4n k_0$. It follows that $n=k^2k_0$.
Now let us compare the inclusions
${\Cal A}_{U(1)_{2k^2k_0}} \subset {\Cal A}$
and ${\Cal A}_{U(1)_{2k^2k_0}} \subset {\Cal A}_{U(1)_{2k_0}}$. Since the 
vacuum representations $H_{\Cal A}$
and $H_{{\Cal A}_{U(1)_{2k_0}}}$ have the same decompositions with respect to
${\Cal A}_{U(1)_{2k^2k_0}}$, we can identify  $H_{\Cal A}$
and $H_{{\Cal A}_{U(1)_{2k_0}}}$ so that we can assume that  ${\Cal A}_{U(1)_{2k^2k_0}}$ is a common conformal subnet of ${\Cal A}$ and  ${\Cal A}_{U(1)_{2k_0}}$
on the same Hilbert space. Now choose unitary 
charged intertwinners $\phi_{i_0}
\in {\Cal A}(I)$ and  $\psi_{i_0}
\in {{\Cal A}_{U(1)_{2k_0}}}(I)$ such that $Ad_{\phi_{i_0}}$ and
$Ad_{\psi_{i_0}}$ induce the same representation $\pi_{i_0}$ of 
 ${\Cal A}_{U(1)_{2k^2k_0}}$. Define a unitary operator $U$ commuting with
${\Cal A}$ such that
$U\phi_{i_0}^m\Omega= \psi_{i_0}^m\Omega, \forall m\in {\Bbb Z}$, where
$\Omega$ is the vacuum vector. One checks easily that
$U{\Cal A}(I)U^*= {{\Cal A}_{U(1)_{2k_0}}}(I), U\Omega=\Omega$. Since
$U$ commutes with ${\Cal A}$, and so it commutes with the action of $PSL(2,{\Bbb R})$, we have  
 $U{\Cal A}(J)U^*= {{\Cal A}_{U(1)_{2k_0}}}(J), \forall J\in {\Cal I}, 
U\Omega=\Omega$, thus proving that $ {\Cal A}\simeq {{\Cal A}_{U(1)_{2k_0}}}$.
\par
The second case when 
${\Cal A}_{c=1}\subset {\Cal A}_{SU(2)_1}^{U(1)} \subset {\Cal A}$ is 
similar as above. 
${\Cal A}_{SU(2)_1}^{U(1)}$ has a continuous
series of irreducible representations labeled by a real number $q$ which
generates a fusion ring that is isomorphic to ${\Bbb R}$. 
By Lemma 4.3
and Cor. 3.9 we have
$H_{{\Cal A}}= \oplus_{q\in {S}} H_q$ where $S\subset {\Bbb R}$ is
an abelian subgroup, and $q^2\in 2{\Bbb Z}, \forall q\in S$. 
If $S=\{0 \}$, then 
${\Cal A}= {\Cal A}_{SU(2)_1}^{U(1)}$. Assume that  $S\neq\{0 \}$, and let $q_0>0$  
be the least positive number in the discrete set $S$. 
It follows that $S=\{{\Bbb Z}q_0 \}$. Let $n$ be the positive integer
such that $ q_0= \sqrt{2n}$. So we have the decompositions
$H_{{\Cal A}}= \oplus_{k\in {\Bbb Z}} F_{ \sqrt{2n}k}$. Compare this with
${\Cal A}_{SU(2)_1}^{U(1)} \subset {\Cal A}_{U(1)_{2n}}$, and a similar argument
using unitary charged intertwinners as above shows that
${\Cal A}\simeq  {\Cal A}_{U(1)_{2n}}$. 
\subheading{4.2.4
${\Cal A}_{c=1}\subset {\Cal A}_{SU(2)_1}^{D_\infty} \subset {\Cal A}$
case}
Let $O$ be the the nontrivial one-dimensional representation 
of $D_\infty$ (the infinite Dihedral group). 
We will also use  $O$  to denote the
corresponding irreducible   representation of 
${\Cal A}_{SU(2)_1}^{D_\infty}$. We note that the conformal dimensions of 
$O$ are integers , and  by [R2], we can 
choose a representative of $[O]$
such that $O^2=id$. Also the braiding operator $\epsilon(O,O)$ is
a scalar with property $\epsilon(O,O)^2=1$. \par 
Let us first consider the case when $\alpha_{O}$ is
localized on $I$, i.e., $\alpha_{O}$ is a DHR representation of 
${\Cal A}$. Apply Doplicher-Roberts reconstruction to ${\Cal A}$ and
 $\alpha_{O}$ as in Prop. 3.8 of [Mu], we get a conformal net
${\Cal A}_1$ such that ${\Cal A}\subset {\Cal A}_1$ and 
${\Cal A}_{SU(2)_1}^{U(1)} \subset   {\Cal A}_1$. By \S4.2.3 
we can identify ${\Cal A}_1$ as ${\Cal A}_{SU(2)_1}^{U(1)}$ or
${\Cal A}_{U(1)_{2n}}, n\in {\Bbb N}$.  If ${\Cal A}_1= {\Cal A}_{SU(2)_1}^{U(1)}$,
we have ${\Cal A}_{SU(2)_1}^{D_\infty}\subset {\Cal A}\subset {\Cal A}_{SU(2)_1}^{U(1)}$
and it follows that $ {\Cal A}\simeq {\Cal A}_{SU(2)_1}^{D_\infty}$ or 
 ${\Cal A}\simeq {\Cal A}_{SU(2)_1}^{U(1)}$ since  ${\Cal A}_{SU(2)_1}^{D_\infty}$ is the
fixed point subnet  under the action of ${\Bbb Z}_2$. If
 ${\Cal A}_1= {\Cal A}_{U(1)_{2n}}$, note that  ${\Cal A}_{SU(2)_1}^{D_\infty}$ is the fixed point subnet  of  ${\Cal A}_{U(1)_{2n}}$
under the action of $D_\infty$, and so ${\Cal A}$ is
the fixed point net under the action of a closed subgroup of
$D_\infty$ as in \S4.2.2. It follows that  ${\Cal A}$ is isomorphic to
either  ${\Cal A}_{SU(2)_1}^{D_\infty},  {\Cal A}_{SU(2)_1}^{U(1)}$ or 
${\Cal A}_{U(1)_{2m}}$ for some $m\in {\Bbb N}$. \par

Let us now return to the general case. Consider the following inclusions
$$ {\Cal A}_{SU(2)_1}^{D_\infty}(I)\subset  {\Cal A}(I) \subset
 {\Cal A}_{SU(2)_1}^{D_\infty}(I')'$$ on $H_{\Cal A}$. 
Let $\gamma_1:  {\Cal A}_{SU(2)_1}^{D_\infty }(I')'\rightarrow  {\Cal A}(I)$ 
be the canonical endomorphism such that its restriction to $ {\Cal A}(I)$
is the canonical endomorphism from $ {\Cal A}(I)$ to 
$ {\Cal A}_{SU(2)_1}^{D_\infty }(I)$
(cf. [LR]). 
Let $v_1\in Hom (id,
O_{\tilde\epsilon} O_{\epsilon})_{{\Cal A}_{SU(2)_1}^{D_\infty}(I)}$ be
a unitary operator, and let $v_O\in {\Cal A}_{SU(2)_1}^{D_\infty}(I')'$ be the unique
unitary operator such that
$ \gamma_1(v_O)= v_1$. Notice that
$\gamma(\tilde\alpha_O\alpha_O)=  O_{\tilde\epsilon} O_{\epsilon}\gamma$,
and it follows that
$v_O a= \tilde\alpha_O\alpha_O(a) v_O, \forall a\in  {\Cal A}(I)$. 
Restrict to $ {\Cal A}_{SU(2)_1}^{D_\infty}(I)$ and recall that $O^2=id$ 
we have that $v_O\in 
{\Cal A}_{SU(2)_1}^{D_\infty}(I)'\cap  {\Cal A}_{SU(2)_1}^{D_\infty}(I')'$. So $v_O$ commutes with
$PSL(2,{\Bbb R})$ by the strong additivity of $ {\Cal A}_{SU(2)_1}^{D_\infty}$. On the
other hand we note that by (2) of  Prop. 3.5 
$\tilde\alpha_O\alpha_O= \alpha_O\tilde\alpha_O$
since $\epsilon(O,O)$ is a scalar, and so $v_O^2 \in {\Cal A}_{SU(2)_1}^{D_\infty}(I)'\cap  {\Cal A}(I)= {\Bbb C} $. 
Multiplying by a scalar if necessary, we can assume that 
$v_O^2=id,  v_O^*=v_O$. Hence the subnet defined by
$\hat{\Cal A}(J):= {\Cal A}(J)\cap \{ v_O \}'$ is a conformal
subnet of ${\Cal A}$ and ${\Cal A}_{SU(2)_1}^{U(1)}\subset \hat{\Cal A}$ is 
a conformal subnet. Now consider 
$\beta_O:= \tilde\alpha_O^{{\Cal A}_{SU(2)_1}^{U(1)}\rightarrow\hat{\Cal A}}
\alpha_O^{{\Cal A}_{SU(2)_1}^{U(1)}\rightarrow\hat{\Cal A}}$. By the definition 
we have 
$v_O   a= \beta_O(a) v_O, \forall a\in  \hat{\Cal A}(I)$,
and by the definition of $\hat{\Cal A}$ we get
$$v_O a= av_O= \beta_O(a) v_O, \forall a\in  \hat{\Cal A}(I)$$
and so $ \beta_O=id$. It follows that  
$$[\tilde\alpha_O^{{\Cal A}_{SU(2)_1}^{U(1)}\rightarrow\hat{\Cal A}}]=
[\alpha_O^{{\Cal A}_{SU(2)_1}^{U(1)}\rightarrow\hat{\Cal A}}]    
$$ 
and so $\alpha_O^{{\Cal A}_{SU(2)_1}^{U(1)}\rightarrow\hat{\Cal A}}$ is localized on $I$. So we can apply the first part of this section to identify
$ \hat{\Cal A}$ as either ${\Cal A}_{SU(2)_1}^{D_\infty},  {\Cal A}_{SU(2)_1}^{U(1)}$  
${\Cal A}_{U(1)_{2m}}^{\Bbb Z_2}$ or ${\Cal A}_{U(1)_{2m}}$ 
for some positive integer $m$. Since $\hat{\Cal A}\subset{\Cal A}$, to
identify $ {\Cal A}$ it is enough to consider the case when 
 $ \hat{\Cal A}\simeq {\Cal A}_{SU(2)_1}^{D_\infty}$ or 
$ \hat{\Cal A}\simeq {\Cal A}_{U(1)_{2m}}^{\Bbb Z_2}$ since the other cases have
been treated in \S4.2.3. Let us show that in this case
$\hat{\Cal A}= {\Cal A}$. Note that by definition
the index of $\hat{\Cal A}(I)\subset  {\Cal A}$ is at most $2$. So we just
have to show that the index is not 2. 
Consider the inclusion $\hat{\Cal A}(I)\subset  {\Cal A}$. 
If the index is 2, one checks easily that
$\alpha_O$ is not localized on $I$ and 
$ [\gamma_{{\Cal A}}]= [1+\tilde\alpha_O\alpha_O]$  

Let $f$ be the vector representation of  $\hat{\Cal A}(I)$. Note that
$[Of]=[f]$. It  is now easy  to check that 
$$
\langle \alpha_f, \alpha_f\rangle = 1, 
\langle \tilde \alpha_f\alpha_f , 1\rangle \geq1
$$
It follows that $ \alpha_f$ is localized on $I$, and 
so $\alpha_O\prec\alpha_f^2 $ is also localized on $I$, a contradiction.  

\subheading{4.2.5
${\Cal A}_{c=1}\subset {\Cal A}_{SU(2)_1}^{G} \subset {\Cal A}, G=E_i, i=6,7,8$
}
Note that since $ {\Cal A}_{SU(2)_1}^{G}$ is absolutely rational, it follows 
that ${\Cal A}$ by is also absolutely rational by Prop. 2.3 of [KL] . If ${\Cal A}\subset {\Cal A}_1$ 
is a conformal subnet of ${\Cal A}_1$ which is irreducible, i.e., 
${\Cal A}(I)'\cap{\Cal A}_1(I)={\Bbb C}$, by Prop. 2.3 of [KL] 
${\Cal A}_1$ is also  absolutely rational, and 
$\mu_{{\Cal A}_1}\leq \mu_{{\Cal A}}$ with equality iff
${\Cal A}= {\Cal A}_1$ by Prop. 2.2. Also note that if
 ${\Cal A}\subset {\Cal A}_1\subset  {\Cal A}_2$ are conformal
subnets with  ${\Cal A}\subset {\Cal A}_1$ and  ${\Cal A}\subset {\Cal A}_2$
irreducible, by  Prop. 2.3 of [KL] again we have that 
 ${\Cal A}\subset {\Cal A}_2$ has finite index and so 
${\Cal A}\subset {\Cal A}_2$ is also irreducible by Lemma 14 of [Ls]. 
We claim that among the irreducible conformal extensions of
${\Cal A}$ there must be a conformal net ${\Cal A}_{max}$ such that
if  ${\Cal A}_{max}\subset  {\Cal B}$ is an irreducible conformal subnet,
then ${\Cal A}_{max}=  {\Cal B}$. If not, we have
${\Cal A}\subset {\Cal A}_1\subset...\subset {\Cal A}_n\subset...
$
where ${\Cal A}_i\subset {\Cal A}_{i+1}$ is irreducible and
${\Cal A}_i\neq {\Cal A}_{i+1}, \forall i$. By Jones's theorem (cf. [J])
the index $[{\Cal A}_{i+1}, {\Cal A}_{i}]\geq 2$, and it follows that
from Prop. 2.2 
$1\leq \mu_{{\Cal A}_{i}} \leq \frac{1}{4^i}\mu_{{\Cal A}}
,\forall i, $
contracting the fact that $\mu_{{\Cal A}}<\infty$. \par
Let ${\Cal A}_{max}$ be a maximal conformal extension of  ${\Cal A}$ 
in the sense above such that
$ {\Cal A}_{SU(2)_1}^{G} \subset {\Cal A}\subset {\Cal A}_{max}$. Let
${\Cal B}$ be the conformal subnet of ${\Cal A}_{max}$ generated by the
$ {\Cal A}_{c=1}$ and 
charged intertwinners
associated with the set $\{ L(1,p^2), \forall p\in {\Bbb Z} \}$. 
 Note that $ {\Cal A}_{SU(2)_1}^{G}\subset {\Cal B}$. 
By \S4.2.2, we can identify ${\Cal B}$ with 
$ {\Cal A}_{SU(2)_1}^{G'}$ where $G'$ a closed A-D-E subgroup of $SO(3)$. 
If $G'$ is of type $D$,  using \S4.2.4   to identify all possible 
${\Cal A}_{max}$, we see that ${\Cal A}_{max}$ can be further
extended, contradicting the maximality of ${\Cal A}_{max}$.\par 
If $G'$ is of type $E$, let $\hat G'$ be the two-fold covering group of
$G'$ in $SU(2)$. For any irreducible representation $\lambda$ of 
$\hat G'$, we use the same $\lambda$ to label the covariant representation of
${\Cal A}_{SU(2)_1}^{G'}$. Note that the statistical dimensions of $\lambda$
are given by the Perron-Frobenius eigenvectors labeled in [GHJ].  
Consider the induction for the pair
${\Cal B}={\Cal A}_{SU(2)_1}^{G'}\subset {\Cal A}_{max}$. By the definition of
${\Cal B}$ and (3) of Lemma 3.3 each $\alpha_\lambda$ is irreducible. 
Consider the index set $S:=\{ \lambda| \lambda\in Irrep G', 
[\alpha_\lambda]= [\tilde \alpha_\lambda] \}$ and 
$\hat S: =\{ \lambda| \lambda\in Irrep \hat G', 
[\alpha_\lambda]= [\tilde \alpha_\lambda] \}$. Note that the set $S$ consists
of  representations of $ {\Cal A}_{max}$ which have  permutation 
symmetry. By Doplicher-Roberts reconstruction as in Prop. 3.9 of [Mu],
we conclude there is an irreducible 
conformal extension of  $ {\Cal A}_{max}$ and it 
follows that $S={1}$ where $1$ denotes the trivial representation 
by the maximality of $ {\Cal A}_{max}$. 
Since $\hat S^2\subset  S$, it follows that there is at most one 
$\lambda\in \hat S-irrep G'$ such that $\lambda^2=1$. By inspecting the
$E$ graph on Page 14 of [GHJ], we conclude that $ \hat S=\{1\}$, and so 

$$
\align
\langle \alpha_\lambda\tilde \alpha_\lambda, 
\alpha_\mu\tilde \alpha_\mu\rangle &=  \langle \alpha_\lambda
\alpha_{\bar \mu},
\tilde \alpha_\mu \tilde\alpha_{\bar \lambda}\rangle \\
&=\langle \sum_{\delta\in irrep \hat G'}
N_{\lambda\bar\mu}^{\delta}\alpha_\delta,
\sum_{\delta'\in irrep \hat G'} N_{\mu\bar\lambda}^{\delta'}\tilde 
\alpha_{\delta'} \rangle\\
&=\delta_{\lambda\mu}
\endalign
$$ 
On the other hand
$\langle \gamma_{\Cal A}, \alpha_\lambda\tilde \alpha_\lambda\rangle
\geq 1$ by Frobenius reciprocity and definitions, and  we have
$$\gamma_{\Cal A}\succ \sum_{\lambda\in irrep \hat G'} 
[\alpha_\lambda\tilde \alpha_\lambda]
$$
This implies
$$
[{\Cal A}_{max}: {\Cal A}_{SU(2)_1}^{G'}]^2 \geq |\hat G'|^2 = 4 |G|^2
$$
and by Prop. 2.2
$$
1\leq \mu_{{\Cal A}_{max}}= \frac{2|G'|^2}
{[{\Cal A}_{max}: {\Cal A}_{SU(2)_1}^{G'}]^2}
\leq \frac{1}{2}
$$
a contradiction. \par
It follows that $G'$ must be of type $A$, and so ${\Cal A}_{max}$ can be
identified with ${\Cal A}_{U(1)_{2m}}$ for some positive integer by \S4.2.3. 
By the maximality of ${\Cal A}_{max}$, either $m=1$ or 
$m$ must be square free, i.e., $m$ is not divisible by $k^2, \forall k>1$. 
In the later case the principal graph of 
$\alpha_{L(1,\frac{1}{4})}^{{\Cal A}\rightarrow {\Cal A}_{max}}$ is
$A_{\infty,-\infty}$ by (3) of Lemma 4.4. On the other hand by 
(3) of Lemma 3.3, $\alpha_{L(1,\frac{1}{4})}^{{\Cal A}\rightarrow {\Cal A}_{max}}$ is also the induced endomorphism of the vector representation $f$ from
$ {\Cal A}_{SU(2)_1}^{G} $ to ${\Cal A}_{max}$, and it has finite depth since
${\Cal A}_{SU(2)_1}^{G}$ is absolutely rational. This contraction
shows that  ${\Cal A}_{max}$ must be identified with 
 ${\Cal A}_{U(1)_{2}}= {\Cal A}_{SU(2)_{1}}$. Now we have inclusions
$$
{\Cal A}_{c=1} \subset {\Cal A}\subset {\Cal A}_{SU(2)_{1}}
$$
Since ${\Cal A}_{c=1}$ is the fixed point subnet of 
${\Cal A}_{SU(2)_{1}}$ under the action of $SO(3)$, it follows that
${\Cal A}$ is the fixed point  subnet of 
${\Cal A}_{SU(2)_{1}}$ under the action of of a closed subgroup of $SO(3)$
as in the end of \S4.2.2. 
\proclaim{Proof of Theorem 4.6}
\endproclaim
If  ${\Cal A}_{c=1}= {\Cal A}$ there is nothing to prove. 
Let us assume that ${\Cal A}_{c=1}\neq {\Cal A}$.
Let  ${\Cal A}_f \subset {\Cal A}$ be the conformal subnet generated by
charged intertwinners associated with the set 
$\{L(1,p^2) , \forall p\in {\Bbb Z},\}$ as in  (3) of 
Cor. 3.9.  
By the spectrum condition the conformal subnet ${\Cal A}_f $
is larger than  ${\Cal A}_{c=1}$, i.e., 
${\Cal A}_{c=1}\subset{\Cal A}_f$ but $ {\Cal A}_{c=1}\neq {\Cal A}_f$. 
By (3) of Cor. 3.9
${\Cal A}_{c=1}\subset{\Cal A}_f$ is discrete in the sense of \S4.2.2,
and by \S4.2.2
${\Cal A}_f$ can be identified as conformal net containing $ {\Cal A}_{SU(2)_1}^{U(1)}$
$ {\Cal A}_{SU(2)_1}^{D_\infty}$ or $ {\Cal A}_{SU(2)_1}^G, G=E_6,E_7,E_8$, and the 
theorem follows from \S4.2.3-5.  \qed 
\par
We note that if the spectrum condition is violated, then the
canonical endomorphism $\gamma$ satisfies strange properties, and
we have been able to rule out some cases in [Xt]. In general it is
still an open question to show that the spectrum condition 
is always satisfied. \par
We point out one application of Th. 4.6. By a general result of  
[Xc] 
the net associated with the coset $SU(2)_4 \subset SU(2)_2\times SU(2)_2$ 
has $13$ irreducible representations, and this net has central 
charge $c=1$. Using the known character formulas in this case one can 
verify that the spectrum condition is
satisfied in this case, so this net must be identified with an element
in the list (*) by Th. 4.6. From the fusion rules of the coset  
in [Xc] one immediately identifies the coset with
${\Cal A}_{U(1)_{12}}^{\Bbb Z_2}$. This observation was first pointed out to me
by C. Dong during our discussions on the ``fixed point resolution''
problem about the cosets, of which   $SU(2)_4 \subset SU(2)_2\times SU(2)_2$
is the first nontrival example (cf. [Xc]).
\heading References \endheading 
\roster 
\item"{[ALR]}"
C. D'Antoni, R. Longo and  F. R\u{a}dulescu, {\it 
Conformal nets, maximal temperature and models from free probability,}  
J. Operator Theory  45  (2001),  no. 1, 195--208. 
\item"[BE1]" J. B\"{o}ckenhauer, D. E. Evans,
{\it Modular invariants, graphs and $\alpha$-induction for
nets of subfactors. I.},
Comm.Math.Phys., {\bf 197}, 361-386, 1998
\item"[BE2]" J. B\"{o}ckenhauer, D. E. Evans,
{\it Modular invariants, graphs and $\alpha$-induction for
nets of subfactors. II.},
Comm.Math.Phys., {\bf 200}, 57-103, 1999   
\item"[BE3]" J. B\"{o}ckenhauer, D. E. Evans,
{\it Modular invariants, graphs and $\alpha$-induction for
nets of subfactors. III.},
Comm.Math.Phys., {\bf 205}, 183-228, 1999
\item"[BMT]" D. Buchholz, G. Mack and I. Todorov, {\it The current 
algebra on the circle as a germ of local field theories}, 
Nucl. Phys. B, Proc. Suppl. {\bf 56}, 20 (1988)
\item"[BS]" D. Buchholz and H. Schulz-Mirbach, {\it Haag duality in conformal
quantum field theory}, Rev. Math. Phys. {\bf 2},105 (1990)
\item"[C]" S. Carpi, {\it Classification of subsystems for Haag-Kastler
Nets generated by c=1 chiral current algebra}, 
Lett. Math. Phys., {\bf 47}, 353-364, 1999
\item"[DVVV]" R. Dijkgraaf, C. Vafa, E. Verlinde and H. Verlinde,
{\it The operator algebra of orbifold models,} 
Comm.Math.Phys., {\bf 123}, 429-487, 1989.
\item"[D]" C. Dong, R. Griess Jr.,{\it 
Rank one lattice type vertex operator algebras and their automorphism
groups},  J. Algebra, {\bf 208}, no. 1, 262-275 (1998)
\item"[DR1]" S. Doplicher and J. E. Roberts, {\it Endomorphisms of $C^*$-algebras, cross products and duality for compact groups,} Ann. Math. {\bf 130}, 75-119 (1989)
\item"[DR2]" S. Doplicher and J. E. Roberts, {\it 
Why there is a field algebra with a compact gauge group describing the superselection structure in particle physics,} Comm.Math.Phys., {\bf 131}, 51-107 (1990).
\item"{[FJ]}" K. Fredenhagen and M. J\"{o}r$\beta$, {\it 
Conformal Haag-Kastler nets, pointlike localized fields and the
existence of opeator product expansions},\par
Comm.Math.Phys., {\bf 176}, 541-554 (1996) 
\item"{[FG]}" J. Fr\"{o}hlich and F. Gabbiani, {\it Operator algebras and
Conformal field theory, } Comm. Math. Phys., {\bf 155}, 569-640 (1993). 
\item"{[FRS]}" J\"{u}rgen Fuchs, Ingo Runkel, Christoph Schweigert, 
{\it 
TFT construction of RCFT correlators I: Partition functions
,} 
hep-th/0204148. 

\item"{[GHJ]}"  F.M.Goodman, P.de la Harpe and V.Jones,
{\it Towers of algebras and Coxeter graphs},
MSRI publications, no.14, 1989

\item"{[GW]}" R. Goodman and N.R. Wallach, {\it Structure and unitary cocycle representations of loop groups and the group of diffeomorphisms of the circle}, 
J. Reine Angew. Math. {\bf 347}, 6910133 (1984)
\item"{[GL1]}"  D.Guido and R.Longo, {\it  The Conformal Spin and
Statistics Theorem},  \par
Comm.Math.Phys., {\bf 181}, 11-35 (1996).
\item"{[GL2]}"  D.Guido and R. Longo,{\it Relativistic invariance
and charge conjugation in quantum field theory},
Comm.Math.Phys., {\bf 148}, 521-551 (1992).
\item"{[H]}" U. Haagerup, {\it On the dual weights for crossed products of 
von Neumann algebras II. Applications of operator valued weights.}  
Math. Scan. 43, 119-140 (1978)
\item"{[I]}" M. Izumi, R. Longo and S. Popa, {\it 
A Galois correspondence for compact groups of automorphisms of
von Neumann Algebras with a generalization to Kac algebras.}
J. Funct. Analysis, {\bf 155}, 25-63 (1998)
\item"{[J]}"  V.F.R. Jones, {Index for subfactors}, 
Invent. Math. {\bf 72}, 1-25 (1983). 
\item"{[K]}" Y. Kawahigashi, {\it Classification of operator algebraic conformal field theories}, math.OA/0211141.
\item"{[KL]}" Y. Kawahigashi and  R. Longo, {\it Classification of local
conformal nets. Case $c<1$}, preprint UTMS 2002-2. 
\item"{[KLM]}" Y. Kawahigashi, R. Longo and M. M\"{u}ger,
{\it Multi-interval Subfactors and Modularity of Representations in
Conformal Field theory},  Comm. Math. Phys., {\bf 219}, 631-669 (2001).
\item"{[L1]}"  R. Longo, {\it Minimal index and braided subfactors}, J.
Funct. Anal., {\bf 109}, 98-112 (1992).
\item"{[L2]}"  R. Longo, {\it Duality for Hopf algebras and for subfactors},
I, Comm. Math. Phys., {\bf 159}, 133-150 (1994).
\item"{[L3]}"  R. Longo, {\it Index of subfactors and statistics of
quantum fields}, I, Comm. Math. Phys., {\bf 126}, 217-247 (1989.
\item"{[L4]}"  R. Longo, {\it Index of subfactors and statistics of
quantum fields}, II, Comm. Math. Phys., {\bf 130}, 285-309 (1990).
\item"{[Ls]}"  R. Longo, {\it Conformal subnets and intermediate subfactors
},  \par
math. OA/0102196, to appear in Comm. Math. Phys.
\item"{[LR]}"  R. Longo and K.-H. Rehren, {\it Nets of subfactors},
Rev. Math. Phys., {\bf 7}, 567-597 (1995).  
\item"{[LRo]}"  R. Longo and J. E. Roberts, {\it A theory of dimension,}
K-Theory {\bf 11}, 103-159 (1997)
\item"{[MS]}" G. Moore and N. Seiberg, {\it Taming the conformal zoo},
Lett. Phys. B  {\bf 220}, 422-430, (1989).
\item"{[Mu]}" M. M\"{u}eger,{\it  On charged fields with group symmetry and degeneracies of Verlinde's matrix $S$}, 
Ann. Inst. H. Poincare Phys. Theor.  71  (1999),  no. 4, 359--394. 
\item"[PP]" M. Pimsner and S. Popa, {\it Entropy and index for  subfactors 
}, Ann. Sci. Ecole Norm. Sup., {\bf 19} (1986)
\item"[Po1]" S. Popa, {\it Classifications of subfactors and their 
endomorphisms}, CBMS Lecture Notes, {\bf 86} (1995)
\item"[Po2]" S. Popa, {\it Classifications of amenable subfactors of type II 
}, Acta Math.,  {\bf 172} (1994)
\item"[PS]" A. Pressley and G. Segal, {\it Loop Groups,} O.U.P. 1986.
\item"[PZ]"  V.B. Petkova, J.-B. Zuber, 
{\it Boundary conditions in charge conjugate sl(N) WZW theories}, Proceedings of the NATO Advanced Research Workshop on  Statistical Field Theories,  Como, 2001, eds. A. Cappelli and G. Mussardo, Kluwer Academic Publishers, 2002, p. 161-170
\item"[R1]" Karl-Henning Rehren,  H.R. Tuneke
{\it Fusion rules for the continuum sectors of the Virasoro algebra with c=1
} Lett.Math.Phys. 53 (2000) 305-312.
\item"[R2]" Karl-Henning Rehren, {\it Braid group statistics and their
superselection rules} In : The algebraic theory of superselection
sectors. World Scientific 1990
\item"[R3]" Karl-Henning Rehren, {\it
News from the Virasoro algebra} Lett.Math.Phys. 30 (1994) 125-130
\item"[SS]" S. Str\u{a}til\u{a}, {\it Modular theory in operator algebras},
Editura Academier, 1981
\item"{[TL]}" V. Toledano Laredo, {\it Fusion of
Positive Energy Representations of $LSpin_{2n}$}.
Ph.D. dissertation, University of Cambridge, 1997
\item"{[W]}"  A. Wassermann, {\it Operator algebras and Conformal
field theories III},  Invent. Math. Vol. 133, 467-539 (1998)
\item"[Xb]" F.Xu, {\it   New braided endomorphisms from conformal
inclusions, } \par
Comm.Math.Phys. 192 (1998) 349-403.
\item"[Xc]" F.Xu, {\it Algebraic coset conformal field theories II}, \par
 math.OA/9903096, Publ. RIMS, vol.35 (1999), 795-824.
\item"[Xi]" F.Xu, {\it 3-manifold invariants from cosets},
math.GT/9907077.    
\item"[Xj]" F.Xu, {\it Jones-Wassermann subfactors for
Disconnected Intervals}, Comm. Contemp. Math. Vol. 2, No. 3 (2000), 307-347
\item"[Xo]" F.Xu, {\it Algebraic orbifold conformal field theories},
Fields Institute Communications, Vol. 30, 429-448, (2001) 
\item"[Xt]" F.Xu, to appear. 
\endroster

\enddocument